\documentclass[english,a4paper,12pt]{article}
\usepackage{amssymb,amsfonts,amsmath,amsthm,graphicx}
\usepackage{verbatim}
\usepackage{multirow}
\usepackage{fullpage}

\usepackage{amssymb,amsfonts,amsmath,amsthm,graphicx}
\usepackage{verbatim}
\usepackage{multirow}
\usepackage{array}
\usepackage{psfrag}
\usepackage[sub,ovp]{psfragx}
\usepackage{tikz}
\usepackage{algorithm}
\usepackage{algorithmic}
\usepackage{tablefootnote}

\usepackage{enumitem}
\usepackage{subfigure}
\setlist{nolistsep}

\def\NN{\mathbb{N}}
\def\RR{\mathbb{R}}

\def\bfm#1{\boldsymbol{#1}}

\makeatletter

\newcommand{\Rmnum}[1]{\expandafter\@slowromancap\romannumeral #1@}
\makeatother
\newcommand{\uu}{\alpha}
\newcommand{\vv}{\beta}
\newcommand{\ww}{\gamma}

\definecolor{purple}{RGB}{180,0,180}
\definecolor{purpleDark}{RGB}{110,0,200}

\newcommand{\corr}[1]{#1}
\newcommand{\corrP}[1]{#1}
\newcommand{\corrPD}[1]{#1}

\newtheorem{thm}{Theorem}
\newtheorem{lem}[thm]{Lemma}

\newtheorem{prop}[thm]{Proposition}

\usepackage{fancyhdr}

\fancypagestyle{specialfooter}{%
  \fancyhf{}
  
  \fancyfoot[R]{July 4, 2017}
}

\begin{document}

\thispagestyle{specialfooter}
\begin{center}
{\LARGE \bf Hermite parametric surface interpolation based on
Argyris element}
\end{center}

\begin{center}
\large Ga\v{s}per Jakli\v{c}$\phantom{}^{\rm (a)}$,
Tadej Kandu\v{c}$\phantom{}^{\rm (b)}$
\end{center}

\small
(a) FGG and IMFM, University of Ljubljana and IAM, University of Primorska, Slovenia\\
(b) INdAM, Unit\`a di Ricerca di Firenze c/o DiMaI ``U. Dini", University of Florence, Italy

\begin{abstract}
In this paper, Hermite interpolation by parametric spline surfaces on triangulations is considered. The splines interpolate points, the corresponding tangent planes and normal curvature forms at domain vertices and approximate tangent planes at midpoints of domain edges. Two variations of the scheme are studied: $C^1$ quintic and $G^1$ octic. The latter is of higher polynomial degree but can approximate surfaces of arbitrary topology.
The construction of the approximant is local and fast. Some numerical examples of surface approximation are presented.
\end{abstract}

\section{Introduction}

In Computer-Aided Geometric Design, one of the fundamental problems is \corr{to construct} 
a parametric spline surface that interpolates prescribed spatial data. The data are usually geometric: points, tangent planes, curvature forms, etc. 
When considering Lagrange geometric interpolation problem,
not much is known on existence or construction of interpolation surfaces \cite{Jaklic-Kozak-Krajnc-Vitrih-Zagar-2006-patch}.

A standard approach 
is to impose $G^1$ smoothness conditions between adjacent triangular patches (see \cite{Smooth3D, Tong-Kim-G1-triang, Hahmann-Bonneau-G1-triang, Farin-02-CAGD} and references therein). Most interpolating schemes of this type are local and can form surfaces of arbitrary topology. One of the main concerns in spline surface construction is how to satisfy nonlinear geometric continuity conditions. The complexity of the problem increases at interior vertices where the smoothness conditions interlace (the vertex enclosure/the twist compatibility problem).
Algorithms usually consist of two steps: construction of a wireframe of interpolation boundary curves and computation of interior control points of the patches \cite{Smooth3D, Hahmann-Bonneau-G1-triang, G1_triang_Bezier, Jaklic-Kanduc-WillmoreCubic}. The schemes are generally fairly complex, usually involving additional subdivision processes \corrP{(each macro patch consists of a few micro patches)}, degree raising or blending techniques. 
In \cite{survey-param-triang}, it was pointed out that many algorithms produce surfaces with unpleasing shapes, e.g., with poor curvature distribution or shape defects. Undesirable shapes are often a result of inappropriate boundary curves of the patches.

One of the most well known and relatively simple $G^1$ interpolation schemes on triangular patches was introduced by Shirman and S\'equin \cite{ShirmanSequin, ShirmanSequin_error} and follows a similar procedure as the one introduced by Farin \cite{Smooth3D}. The method interpolates points and tangent planes at the vertices, and consists of quartic patches on Clough--Tocher split. A method by Hahmann and Bonneau solves the vertex enclosure problem by introducing the so-called 4-split \cite{Hahmann-Bonneau-G1-triang}. Although the construction is focused on obtaining good approximation surface, it is not clear how to properly set shape parameters and the number of control points is relatively big considering that the scheme interpolates only points at triangle vertices. In \cite{Tong-Kim-G1-triang}, the authors \corrP{Tong and Kim} consider 
interpolation of points, tangent planes and normal curvatures. 
However, \corrP{they presume} that the approximated surface is given in the implicit form and so additional \corrP{approximation points} are extracted \corrP{and used in a least squares data fitting}. 

An alternative to geometric continuity is to construct splines satisfying stricter $C^r$ continuity conditions \cite{Schumaker-Fasshauer-ParamSpline, Farin-02-CAGD, ME_min_energy, Baramidze-MinE_spherical_splines}. 
The advantage of this approach is that the smoothness conditions are linear and they imply a simple geometric construction of control points.
The main drawbacks are that the schemes cannot approximate a surface of arbitrary topology \cite{Herron-closedSurf} and that for the most interesting
low degrees the dimension of the spline space is still unknown \cite{Schumaker-Lai-splines-triangulations-07,Jaklic-dimension-05}.

Macro-elements are a special type of $C^r$ smooth interpolation splines, defined on triangulated domains  \cite{Schumaker-Lai-splines-triangulations-07, Schumaker-Lai-CT, Schumaker-Alfeld-PS, Clough-Tocher_FE}. Their structure overcomes the problems with the spline space dimension. 
Furthermore, the shape of the spline depends only on local data. The approximants are obtained in a closed form and have the optimal approximation order. 

In this paper, we present an interpolation scheme for parametric surfaces that is based on the $C^1$ polynomial macro-element, known as (quintic) Argyris element \cite{Ciarlet_FEM_1978, Zlamal_FE, MorganScott_NodalBasis_1975, Schumaker-Lai-splines-triangulations-07}. Two variants of the scheme are derived: $C^1$ quintic and $G^1$ octic scheme.
The approximants interpolate given geometric data: points, tangent planes and normal curvature forms. 
The interpolation conditions do not fully determine the shape of the spline. Thus an approach for computing appropriate free shape parameters is introduced.

As a first step of the scheme, a referential linear interpolating surface is constructed. To improve the quality of the surface, one step of the improved Butterfly scheme that can handle arbitrary topology is applied on the control points of the linear spline \cite{Butterfly_improved_Zorin, Dyn_butterfly_1990}.


In order to satisfy interpolation conditions, referential control points are projected onto the corresponding tangent planes. To overcome the twist compatibility problem when enforcing smoothness conditions between the patches, $C^2$ smoothness conditions are imposed at every patch vertex. Corrections of the control points are computed as the solution of a \corrP{small} least squares minimization problem that enforces $C^2$ smoothness.


The construction of the interpolants is local. 
The wireframe of boundary curves is constructed using 
referential control points that better represent \corr{the basic} shape characteristics of the resulting surface. 
Higher polynomial degrees are needed to satisfy \corrP{the} smoothness conditions. The parametric scheme requires degree $8$ to enforce $G^1$ cross-boundary smoothness.  \corr{The polynomial degree can be reduced to 7 \corrP{(or lower)} if certain geometric conditions on the interpolation data are satisfied}.

The paper is organized as follows. In Section~2, basic notation is introduced. 
$C^r$ continuity conditions across a common edge of adjacent patches and at a vertex are recalled in Section~3. In Section~4, $G^1$ smoothness conditions across edges are examined in detail. 
Geometric conditions for reducing the polynomial degree 8 are derived. Construction of control points imposed by three types of interpolation
conditions is analyzed in Section~5. The construction is split into: interpolation of tangent planes at the vertices, interpolation of normal curvature forms at the vertices and approximation of tangent planes at edge midpoints. 
In Section~6, some numerical examples of surface approximation are presented. At the end, main conclusions are emphasized.

\section{Notation}

Let $\triangle$ be a triangulation of a given domain $\Omega \subset \RR^2$. Every edge $e$ and triangle $\tau$ of $\triangle$ is described as a list of vertices $v$: $e=(v_0, v_1)$ and $\tau=(v_0,v_1,v_2)$, respectively. Let the set of all vertices 
be denoted by $\mathcal V$. 
%
In our $G^1$ approximation scheme we will construct only local domain triangulations around interpolation points in order to apply $C^2$ smoothness conditions at the vertices.



Let $\tau\in\triangle$ be a non-degenerate triangle. Every point $v\in \RR^2$ can be written in barycentric coordinates with respect to $\tau$ as $v := v(\tau) := (\uu,\vv,\ww)$, $\uu+\vv+\ww=1$. The Bernstein basis polynomials of total degree $d$ 
are defined as
\begin{align*}
B_{\bfm i}^d(v):=B_{ijk}^d(\uu,\vv,\ww):=
\frac{d!}{i!j!k!} \, \uu^i\vv^j\ww^k,\qquad |\bfm i|=d.
\end{align*}

A parametric polynomial $\bfm p$ of total degree $d$ 
can be represented in the B\'{e}zier form
\begin{align*}
\bfm p=\sum_{|\bfm i|=d}\bfm c_{\bfm i} B_{\bfm i}^d,
\end{align*}
where $\bfm c_{\bfm i}=\bfm c_{ijk}\in\RR^3$ are its control points. 

\begin{figure}[!htb]
\centering
\begin{overpic}[width=6cm]{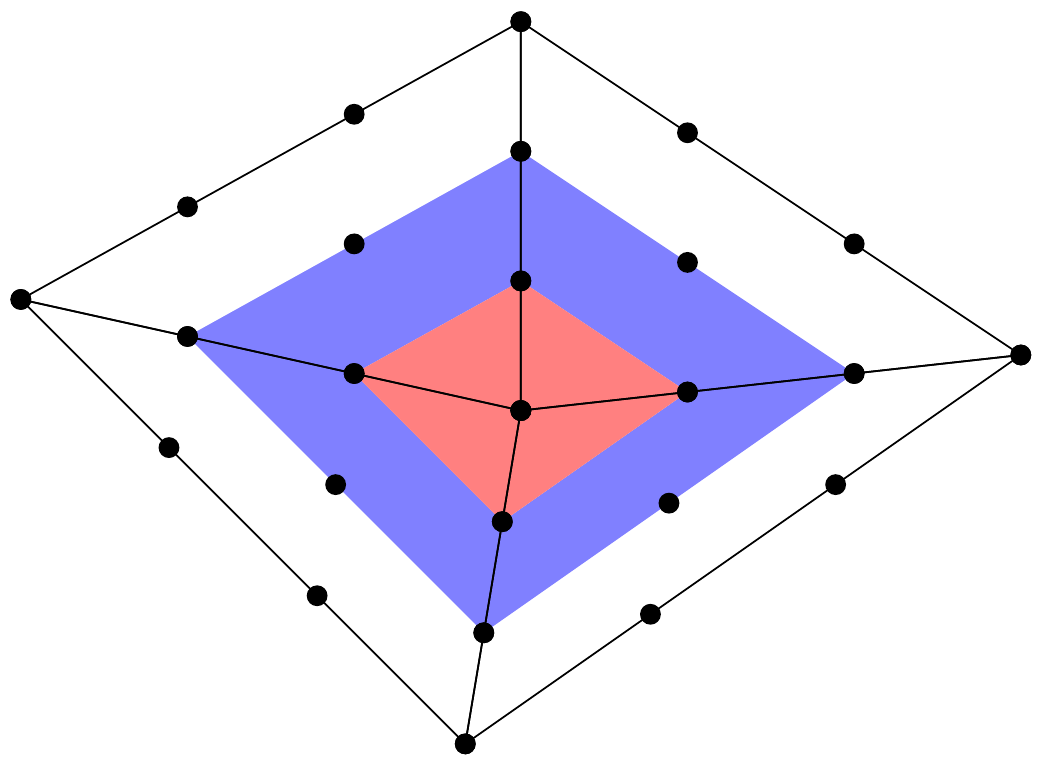}
\put(50.5,29.25){\small $\bfm P$}
\put(41.5,37){\small $\mathcal D_1(v)$}
\put(41.5,49.5){\small $\mathcal D_2(v)$}
\end{overpic}
\caption{Set $\mathcal D_0(v) = \{\bfm P\}$ consists of the control point in the center $v$. Sets $\mathcal D_1(v)$ and $\mathcal D_2(v)$ are represented with black dots in red and red+blue area, respectively. 
}
\label{fig:setDvAndDe}
\end{figure}

\emph{Disk} $\mathcal D_\ell(v)$, $\ell \geq 0$, is a set of control points \corr{of a spline} that are at most $\ell$ indices away from the origin $v\in\mathcal V$ (see Fig.~\ref{fig:setDvAndDe}). \emph{Ring} is defined as $\mathcal R_\ell(v) := \mathcal D_\ell(v) \backslash \mathcal D_{\ell-1}(v)$ for $\ell \geq 1$. We will always presume $C^0$ continuity.

\corr{For a
vector of scalars} $\bfm a= (a_\ell)_{\ell=1}^r$ and \corr{a vector} $\bfm b= (b_\ell)_{\ell=1}^r$,  consisting of scalars or points, \corr{we define a scalar product as} 
\begin{align*}
\langle \bfm a, \bfm b \rangle := \sum_{\ell=1}^r a_\ell \, b_\ell.
\end{align*}

Before constructing our spline interpolant, a referential spline surface that interpolates given data points \corr{is} constructed. In our scheme we presume that a spatial triangulation (i.e., a linear spline interpolant) passing through the interpolation points is already given. After that, one step of the modified Butterfly scheme is applied on control points of the linear spline \cite{Butterfly_improved_Zorin}. 
\corr{That way we obtain a better starting approximation surface that combines} data also from the neighbouring patches. 
The symbol $\bullet^{\bowtie}$ will be used to indicate different objects (patches, control points, sets) that correspond to the referential interpolant. Polynomial degree of the obtained quadratic patches \corr{needs to be raised} to 5 for the $C^1$ and to 8 for the $G^1$ scheme.

\section{$C^r$ smooth splines}

A spline $\bfm s $ 
consists of patches $\bfm p^{[\tau]}$, $\bfm s|_\tau =: \bfm p^{[\tau]}= \sum \bfm c_{\bfm i}^{[\tau]} B_{\bfm i}^d$, for $\tau\in\triangle$.
Let $\bfm e_1,\, \bfm e_2,\, \bfm e_3$ be $(1,0,0)$, $(0,1,0)$ and $(0,0,1)$, respectively.
The intermediate de Casteljau points for parameter $v=(\uu,\vv,\ww)$ are defined as
\begin{align*}
\bfm c_{\bfm i}^{(k)} := \bfm c_{\bfm i}^{(k)}(v) := \left\langle v, \left(\bfm c_{\bfm i+\bfm e_1}^{(k-1)}, \bfm c_{\bfm i+\bfm e_2}^{(k-1)}, \bfm c_{\bfm i+\bfm e_3}^{(k-1)}\right) \right \rangle,\qquad |\bfm i| = d - k,
\end{align*}
and $\bfm c_{\bfm i}^{(0)} := \bfm c_{\bfm i}$.

The following two well known theorems state the $C^r$ continuity conditions across an adjoining edge and at a vertex \cite{Schumaker-Lai-splines-triangulations-07, Farin-02-CAGD}.

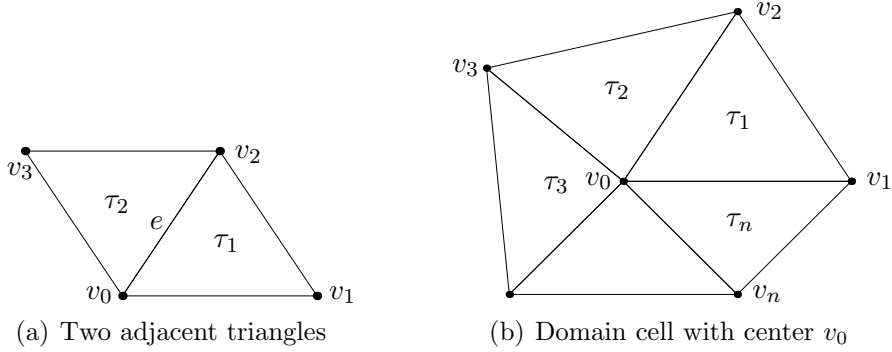
\begin{figure}[!htb]
\centering
\subfigure[Two adjacent triangles]{
\begin{overpic}[width=4cm]{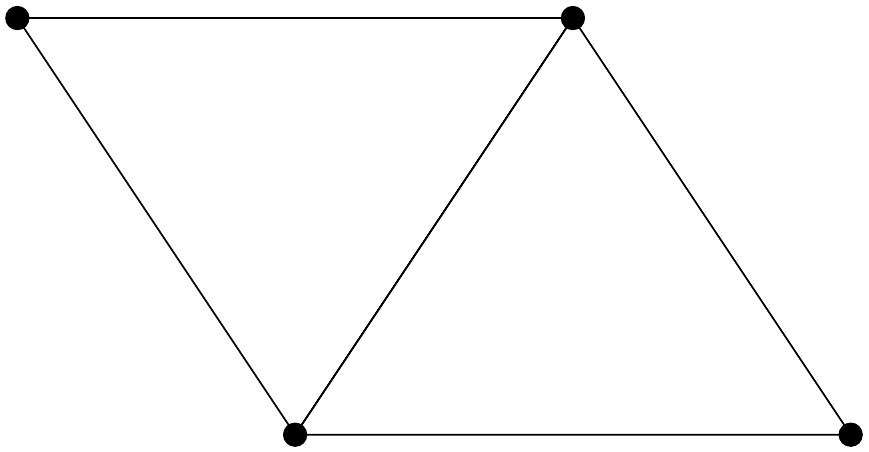}
\put(22,1){\small $v_0$}
\put(102,1){\small $v_1$}
\put(71,47.5){\small $v_2$}
\put(-4,42){\small $v_3$}
\put(43,24){\small $e$}
\put(64,18){\small $\tau_1$}
\put(28,31){\small $\tau_2$}
\end{overpic}
}
\hspace{1.5cm}
\subfigure[\corrP{Domain} cell with center $v_0$]{
\begin{overpic}[width=5cm]{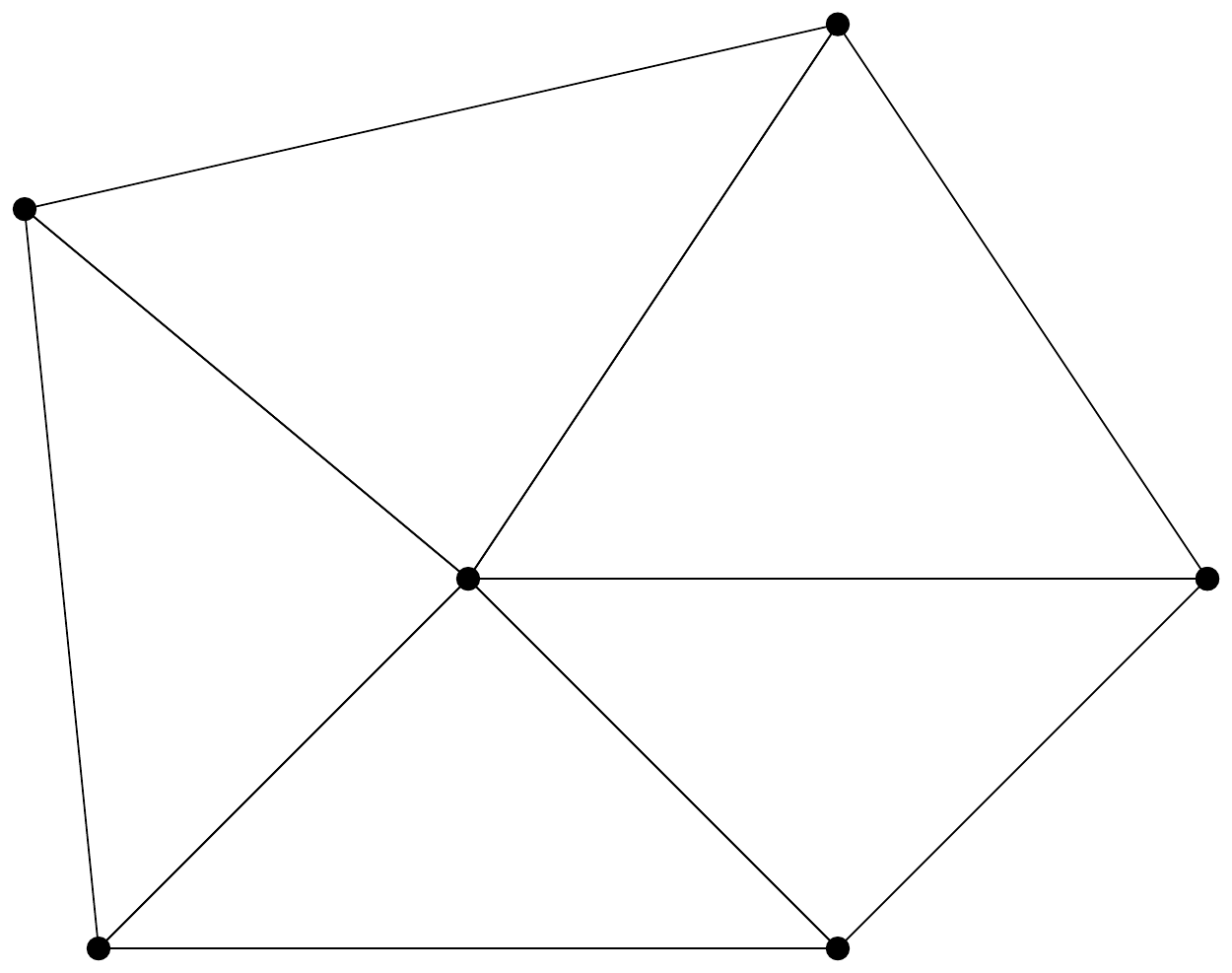}
\put(28,31){\small $v_0$}
\put(102,31){\small $v_1$}
\put(65,47){\small $\tau_1$}
\put(73,76){\small $v_2$}
\put(33,55){\small $\tau_2$}
\put(-7,61){\small $v_3$}
\put(17,30){\small $\tau_3$}
\put(72,1){\small $v_n$}
\put(65,20){\small $\tau_n$}
\end{overpic}
}
\caption{Two sets of adjacent triangles relevant for smoothness conditions.}
\label{fig:triangs}
\end{figure}

\begin{thm}\label{thm:smoothEdge}
Let $\bfm p^{[\tau_1]}$ 
and $\bfm p^{[\tau_2]}$ 
be adjacent patches, defined on triangle $\tau_1=(v_0,v_1,v_2)$ and $\tau_2=(v_0,v_2,v_3)$, respectively (see Fig.~\ref{fig:triangs}(a)). For $0\leq r\leq d$, the patches join with $C^r$ continuity across the edge $e=(v_0,v_2)$ if
\begin{align*}
\bfm c_{ijk}^{[\tau_2]}=\left(\bfm c_{i0j}^{[\tau_1]}\right)^{(k)}(v_3(\tau_1)),\qquad k \leq r,\;|\bfm i|=d.
\end{align*}
\end{thm}

\begin{thm}\label{thm:smoothPoint}
Let $\triangle$ be a 
triangulation with triangles $\{\tau_{\ell}=( v_0,v_\ell,v_{\ell+1})\}_{\ell=1}^{n}$ (Fig.~\ref{fig:triangs}(b)). 
If $v_0$ is interior vertex, $v_{n+1}\equiv v_1$. 
For $0\leq r\leq d$, the patches $\bfm p^{[\tau_\ell]}$ join with $C^r$ continuity at the vertex $v_0$ if
\begin{align*}
\bfm c_{ijk}^{[\tau_{\ell+1}]}=\left(\bfm c^{[\tau_\ell]}_{i0j}\right)^{(k)}(v_{\ell+2}(\tau_\ell)),\qquad 
j+k\leq r,\;\; |\bfm i|=d,\;\; 1 \leq \ell\leq n-1.
\end{align*}
\end{thm}
We call a set of triangles in 
\corr{Fig.~\ref{fig:triangs}(b)} 
a \emph{\corrP{domain} cell}.
%
%
%

\section{$G^1$ geometric smoothness}

%
To construct $G^1$ octic interpolant, first we need to analyze geometric smoothness conditions across common boundary curves of the adjacent patch\-es.
Let $\tau_1$ and $\tau_2$ be adjacent triangles as in \corr{Fig.~\ref{fig:triangs}(a)} 
and let $v(t)=(1-t)\,v_0 + t\,v_2$ for $t\in[0,1]$. Let $\bfm b(t) := D_{v_2-v_0} \bfm p^{[\tau_1]}(v(t))$ be a directional derivative of $\bfm p^{[\tau_1]}$ along the common boundary curve 
\corr{and let $\bfm v$ be an unknown \emph{transversal vector function}.}
The patches $\bfm p^{[\tau_1]}$ and $\bfm p^{[\tau_2]}$ join with $G^1$ geometric continuity if there exist \emph{connecting functions} $\lambda, \mu, \nu, \xi$ \corr{and $\bfm v$} such that 
\begin{align}\label{eqn:G1_Chiy}
\nonumber D_{v_3-v_0} \bfm p^{[\tau_2]}(v) &= \lambda(t)\bfm b(t)  + \mu(t) \bfm v(t),\\[-1.6ex]
&& t \in [0,1].\\[-1.6ex]
\nonumber D_{v_1-v_0} \bfm p^{[\tau_1]}(v) &= \nu(t)\bfm b(t) + \xi(t) \bfm v(t),
\end{align}
The transversal vector function $\bfm v$ and the boundary vector function $\bfm b$ span the tangent plane of $\bfm p^{[\tau_1]}$ and $\bfm p^{[\tau_2]}$ at vertex $v$. In practice, the connecting functions \corr{and $\bfm v$} are \corr{(parametric)} polynomials of prescribed degree. If all of the connecting functions $\lambda, \mu, \nu, \xi$ are constant we obtain $C^1$ smoothness conditions. 

\begin{figure}[!htb]
\centering
\begin{overpic}[trim = 0cm 3.75cm 0cm 3.75cm, clip=true, width=6cm]{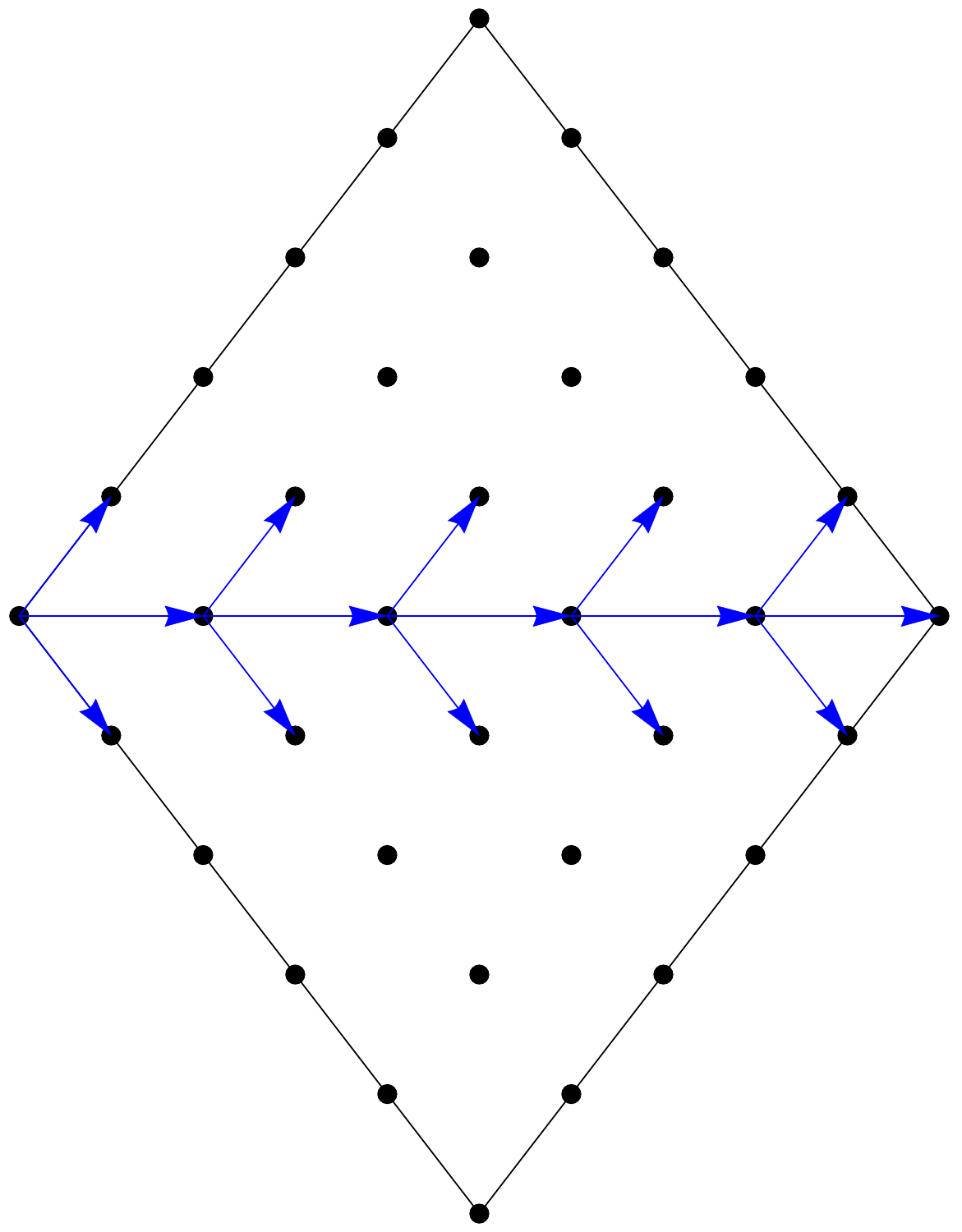}
\put(13,35){\small $\bfm b_0$}
\put(32,35){\small $\bfm b_1$}
\put(51,35){\small $\bfm b_2$}
\put(70,35){\small $\bfm b_3$}
\put(89,35){\small $\bfm b_4$}

\put(1,43){\small $\bfm d_0$}
\put(20,43){\small $\bfm d_1$}
\put(39,43){\small $\bfm d_2$}
\put(59,43){\small $\bfm d_3$}
\put(78,43){\small $\bfm d_4$}

\put(1,18){\small $\bfm e_0$}
\put(20,18){\small $\bfm e_1$}
\put(39,18){\small $\bfm e_2$}
\put(59,18){\small $\bfm e_3$}
\put(78,18){\small $\bfm e_4$}
\end{overpic}
\hspace{.25cm}
\raisebox{.45cm}{\begin{overpic}[width=7cm,]{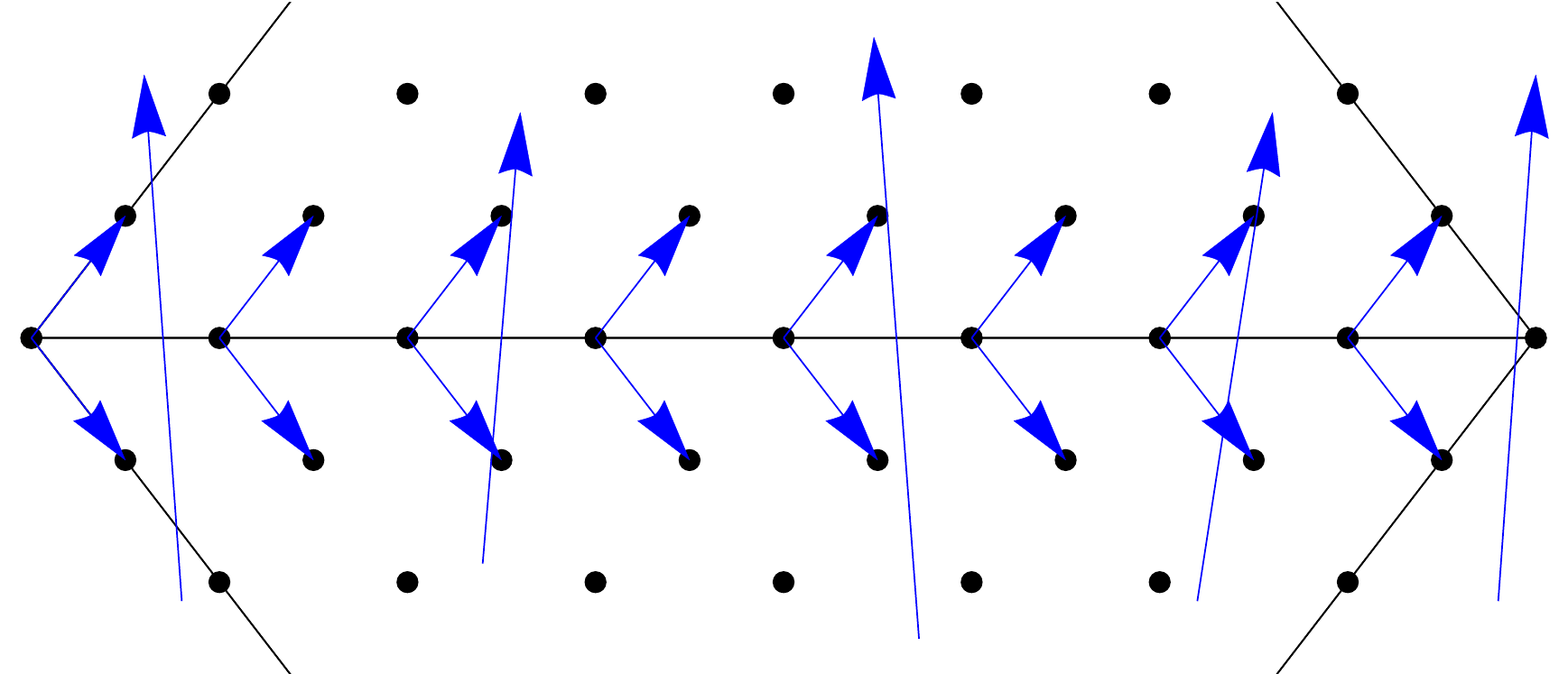}
\put(8,45){\small $\bfm v_0$}
\put(30,42){\small $\bfm v_1$}
\put(50,46){\small $\bfm v_2$}
\put(77,40){\small $\bfm v_3$}
\put(95,42){\small $\bfm v_4$}

\put(-2,26){\small $\bfm d'_0$}
\put(12,26){\small $\bfm d'_1$}
\put(24,26){\small $\bfm d'_2$}
\put(35,26){\small $\bfm d'_3$}
\put(47,26){\small $\bfm d'_4$}
\put(60,26){\small $\bfm d'_5$}
\put(71,26){\small $\bfm d'_6$}
\put(83,26){\small $\bfm d'_7$}

\put(-2,15){\small $\bfm e'_0$}
\put(12,15){\small $\bfm e'_1$}
\put(24,15){\small $\bfm e'_2$}
\put(35,15){\small $\bfm e'_3$}
\put(47,15){\small $\bfm e'_4$}
\put(60,15){\small $\bfm e'_5$}
\put(71,15){\small $\bfm e'_6$}
\put(83,15){\small $\bfm e'_7$}
\end{overpic}}
\caption{Vectors $\bfm d_\ell, \bfm e_\ell, \bfm b_\ell, \bfm v_\ell, \bfm d'_\ell, \bfm e'_\ell$  for two adjoining patches.}
\label{fig:g1_pog}
\end{figure}
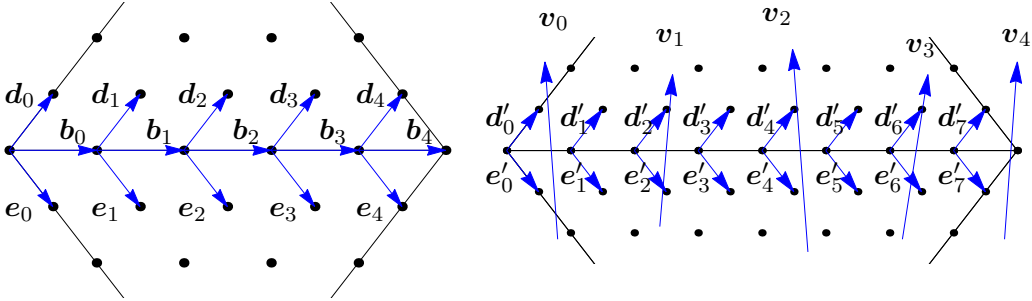

Let the connecting functions be of degree $r$. The functions can be expressed in B\'ezier form using univariate Bernstein polynomials $B_m^r$:
\begin{align*}
\lambda =:\sum_{m=0}^{r} \lambda_m B_{m}^{r},\quad
\mu =:\sum_{m=0}^{r} \mu_m B_{m}^{r},\quad
\nu =:\sum_{m=0}^{r} \nu_m B_{m}^{r},\quad
\xi =:\sum_{m=0}^{r} \xi_m B_{m}^{r}.
\end{align*}
Similarly, let us express the parametric polynomials:
\begin{align*}
D_{v_3-v_0} \bfm p^{[\tau_2]} &=:\sum_{\ell=0}^{d-1} \bfm d_\ell B_{\ell}^{d-1} = \sum_{\ell=0}^{d-1+r} \bfm d'_\ell B_{\ell}^{d-1+r}, &
\bfm b &=:\sum_{\ell=0}^{d-1} \bfm b_\ell B_{\ell}^{d-1},\\
D_{v_1-v_0} \bfm p^{[\tau_1]} &=:\sum_{\ell=0}^{d-1} \bfm e_\ell B_{\ell}^{d-1} = \sum_{\ell=0}^{d-1+r} \bfm e'_\ell B_{\ell}^{d-1+r}, &
\bfm v &=:\sum_{\ell=0}^{d-1} \bfm v_\ell B_{\ell}^{d-1}.
\end{align*}
The control points $\bfm d'_\ell, \bfm e'_\ell$ are obtained from $\bfm d_\ell, \bfm e_\ell$ after raising the polynomial degree $r$ times. 
An example of vectors $\bfm d_\ell, \bfm e_\ell, \bfm b_\ell, \bfm v_\ell, \bfm d'_\ell, \bfm e'_\ell$ for $d=5$ and $r=3$ is depicted in Fig.~\ref{fig:g1_pog}.

By comparing the coefficients on both sides of the system \eqref{eqn:G1_Chiy} we obtain the following $d+r$ vector equations:
\begin{align}
\label{eqn:G1_cond1} \bfm d'_{\ell} &= \sum_{m=0}^r \frac{\binom{d-1}{\ell-m} \binom{r}{m}}{\binom{d-1+r}{\ell}} \left( \lambda_{m}\bfm b_{\ell-m} + \mu_{m}\bfm v_{\ell-m}\right),\\[-1.6ex]
\nonumber & & \ell = 0,1,\dots,d-1+r.\\[-1.6ex]
\label{eqn:G1_cond2} \bfm e'_{\ell} &= \sum_{m=0}^r \frac{\binom{d-1}{\ell-m} \binom{r}{m}}{\binom{d-1+r}{\ell}} \left( \nu_{m}\bfm b_{\ell-m} + \xi_{m}\bfm v_{\ell-m}\right),
\end{align}

To reduce the number of different cases we would need to examine, we presume a common geometrical situation that the vectors $\bfm e_\ell \times \bfm b_\ell$ and $\bfm b_\ell \times \bfm d_\ell$ point in the same direction:
\begin{align}\label{eqn:geom_omejitev}
\frac{1}{\| \bfm e_\ell \times \bfm b_\ell  \|} \bfm e_\ell \times \bfm b_\ell = 
\frac{1}{\| \bfm b_\ell \times \bfm d_\ell \|} \bfm b_\ell \times \bfm d_\ell, \qquad \ell=0,
d-1.
\end{align}
To maintain orientation of the surface, the conditions $\mu>0$ and $\xi < 0$ must be fulfilled.

\subsection{Conditions on the connecting functions} 

From now on let us consider only the case $d=5$. Let us also presume that sets of control points $\mathcal D_2(v_0)$ and $\mathcal D_2(v_2)$ satisfy $C^2$ conditions (Thm.~\ref{thm:smoothPoint}). Therefore, control points $\bfm d'_{\ell},\, \bfm e'_{\ell}$ in \eqref{eqn:G1_cond1} and \eqref{eqn:G1_cond2} are fixed for $\ell=0,1, 3+r, 4+r$. All of the control points $\bfm b_\ell$ are also fixed. Let us show that in order to solve \eqref{eqn:G1_cond1} and \eqref{eqn:G1_cond2} it is sufficient that $\lambda, \mu, \nu, \xi$ are cubic polynomials, i.e., $r=3$. Hence the polynomial degree of the quintic patches is raised to 8.


\begin{thm}
Let $r=3$ and points in $\mathcal D_2(v_0)$ and $\mathcal D_2(v_2)$ be fixed and satisfy $C^2(v_0)$ and \corr{$C^2(v_2)$} \corrP{smoothnesses}. For $\lambda_m \in \RR,\, \mu_m >0$, $m=0,1,2,3$, and $\bfm v_2\in \RR^3$ there exist unique $\nu$, $\xi$ \corr{and $\bfm v$} that satisfy conditions \eqref{eqn:G1_cond1} and \eqref{eqn:G1_cond2}. Control coefficients of the functions $\xi$ and $\nu$ are
\begin{align}\label{eqn:kubicne_povezave}
\nonumber &\xi_0 
= f_0 \, \mu_0,
\nonumber &\nu_0 = 
    f_0 \, \lambda_0 + g_0,\\
\nonumber &\xi_1 = f_0 \, \mu_1,
%
 \nonumber &\nu_1 = f_0 \, \lambda_1 + g_0 ,\\[-1.6ex]
    %
\\[-1.6ex]
\nonumber &\xi_2 = f_4 \, \mu_2,
%
\nonumber &\nu_2 = f_4 \, \lambda_2 + g_4,\\
   %
\nonumber &\xi_3 = f_4 \, \mu_3,
\nonumber &\nu_3 = f_4 \, \lambda_3 + g_4,
\end{align}
with the following parameters
\begin{align}\label{eqn:f0,f4}
\nonumber &f_0 := \frac{\langle \bfm b_0 \times \bfm d_0 , \bfm b_0 \times \bfm e_0 \rangle}{\| \bfm b_0 \times \bfm d_0 \|^2},
&g_0 := -\frac 5 8 \frac{\langle \bfm b_0 \times \bfm d_0 , \bfm d_0 \times \bfm e_0 \rangle}{\| \bfm b_0 \times \bfm d_0 \|^2}, \\[-1.6ex]
\\[-1.6ex]
\nonumber&f_4 := \frac{\langle \bfm b_4 \times \bfm d_4 , \bfm b_4 \times \bfm e_4 \rangle}{\| \bfm b_4 \times \bfm d_4 \|^2},
&g_4 :=  -\frac 5 8 \frac{\langle \bfm b_4 \times \bfm d_4 , \bfm d_4 \times \bfm e_4 \rangle}{\| \bfm b_4 \times \bfm d_4 \|^2}.
\end{align}

\end{thm}

To ensure that the patches lie on the correct side of the half-space, extra conditions $\mu>0$, $\xi<0$ must hold true. The next proposition simplifies the verification of these conditions.

\begin{prop}\label{prop:sign}
$\mu_\ell >0$ if and only if $\xi_\ell < 0$, $\ell = 0,1,2,3$.
\end{prop}


Proposition \ref{prop:sign} simplifies the conditions on $\bfm v$ since it is enough to check the sign of only one out of the two connecting functions $\mu, \xi$.
For example, a simple heuristic way to set the vectors $\bfm v_\ell$:
\begin{align}\label{eqn:simpleVectorV}
\bfm v_\ell := {\bfm d}_\ell - {\bfm e}_\ell, \qquad \ell = 0,1, 3,4,
\end{align}
satisfies conditions \eqref{eqn:G1_cond1} and \eqref{eqn:G1_cond2}. For this case, the conditions on connecting functions simplify considerably and the following relations are obtained
\begin{align*}
\lambda_0 = \lambda_1,\; 
\lambda_2 = \lambda_3, \quad
\mu_0 = \mu_1,\;
\mu_2 = \mu_3, \quad
\nu_0 = \nu_1,\;
\nu_2 = \nu_3, \quad
\xi_0 = \xi_1,\;
\xi_2 = \xi_3.
\end{align*}
The remaining vector $\bfm v_2$ does not affect the smoothness conditions and remains as an additional parameter. It can be used as a shape parameter or to approximate additional data at the interior of the edge. Details of how to set $\bfm v_2$ so that the patches approximate tangent planes at the middle of the edges are explained in Section~5.3. 

\subsection{Reducing the degree of the connecting functions}

Till now the \corr{considered} connecting functions $\lambda, \mu, \nu, \xi$ were cubic polynomials. A natural question arises: Can we reduce the polynomial degree, since we had several free parameters in the cubic case? The answer is in the affirmative and the functions can be quadratic under certain geometric conditions. \corr{This implies that we can employ patches of degree 7 rather than 8 while preserving the same boundary curves and the amount of approximation data}.

%
If $f_0 \neq f_4$, the following \corrP{family of solutions} exists:
\begin{align*}
\mu_1 = \frac{1}{3} \mu_0, \quad 
\mu_2 = \frac{1}{3} \mu_3, \quad
\lambda_1 = \frac{1}{3} \lambda_0 - \frac{2(g_0 - g_4)}{3 (f_0 - f_4)}, \quad 
\lambda_2 = \frac{1}{3} \lambda_3 -  \frac{2(g_0 - g_4)}{3 (f_0 - f_4)}.
\end{align*}
%
When $f_0=f_4$ \corrP{a} solution exists only if $g_0 = g_4$:
\begin{align}\label{eqn:c1_pog}
\mu_0 - 3 \mu_1 + 3 \mu_2 - \mu_3 = 0, \qquad \lambda_0 - 3 \lambda_1 + 3 \lambda_2 - \lambda_3 = 0.
\end{align}
If  $f_0=f_4$ and $g_0 \neq g_4$ both functions $\lambda, \nu$ can not be simultaneously quadratic. When $|g_0-g_4|/|f_0-f_4| \gg 1$ it is better to avoid \corr{using} the solution since big oscillations of the connecting functions can lead to undesired shape defects of the patches.

\corrP{Geometrically, the conditions $f_0=f_4$, $g_0 = g_4$ hold true if there exists an underlying domain triangulation. In this case constant connecting functions satisfy conditions \eqref{eqn:c1_pog} and we get $C^1$ smoothness conditions.}

\section{Interpolation conditions and minimizing \corrPD{rings}}\label{sec:IntApp}

\begin{figure}[!htb]
\centering
\subfigure[$C^1$ quintic patch]{
\begin{overpic}[height=4.5cm]{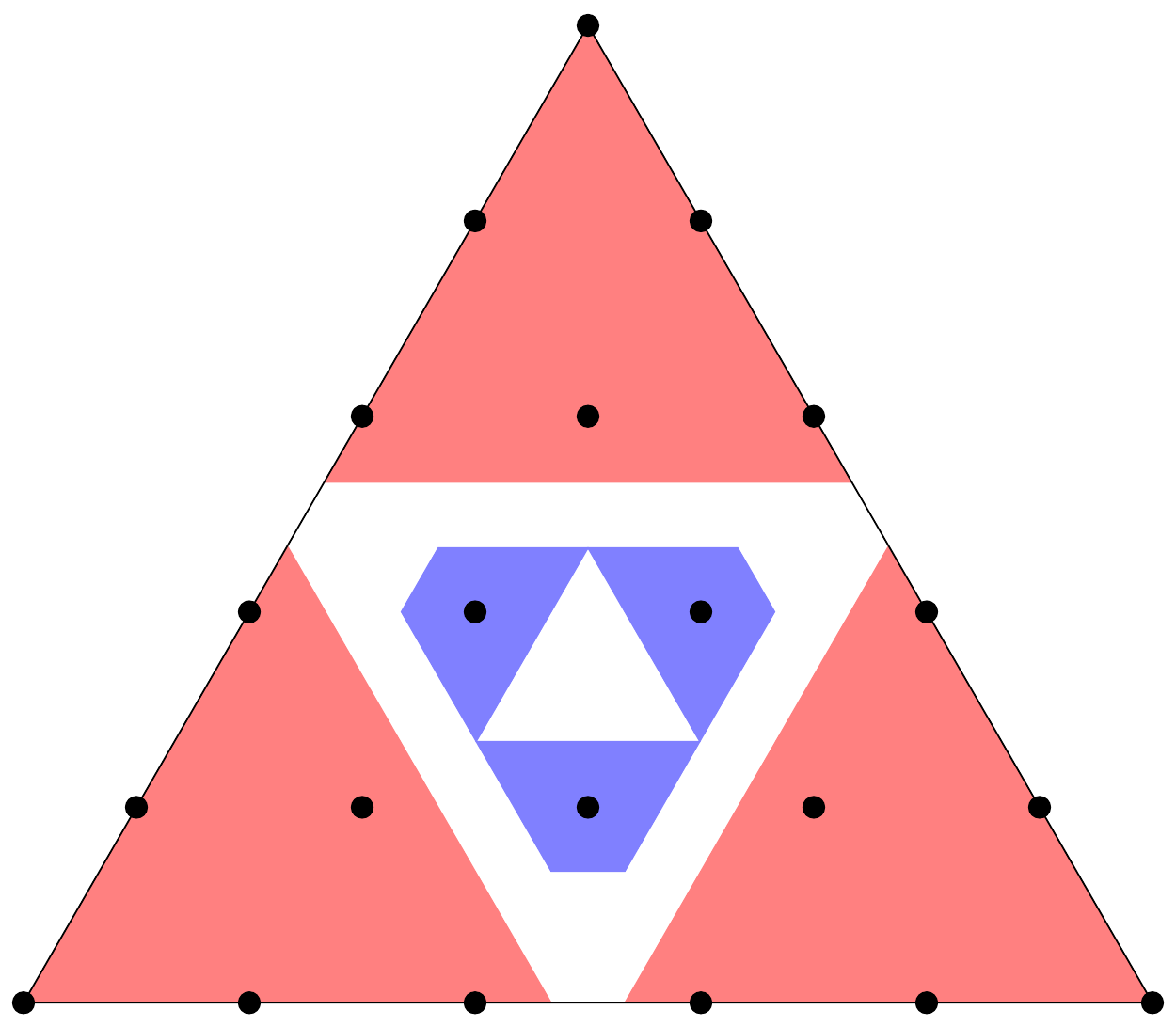}
\end{overpic}
}
\hspace{1cm}
\subfigure[$G^1$ octic patch]{
\begin{overpic}[height=4.5cm]{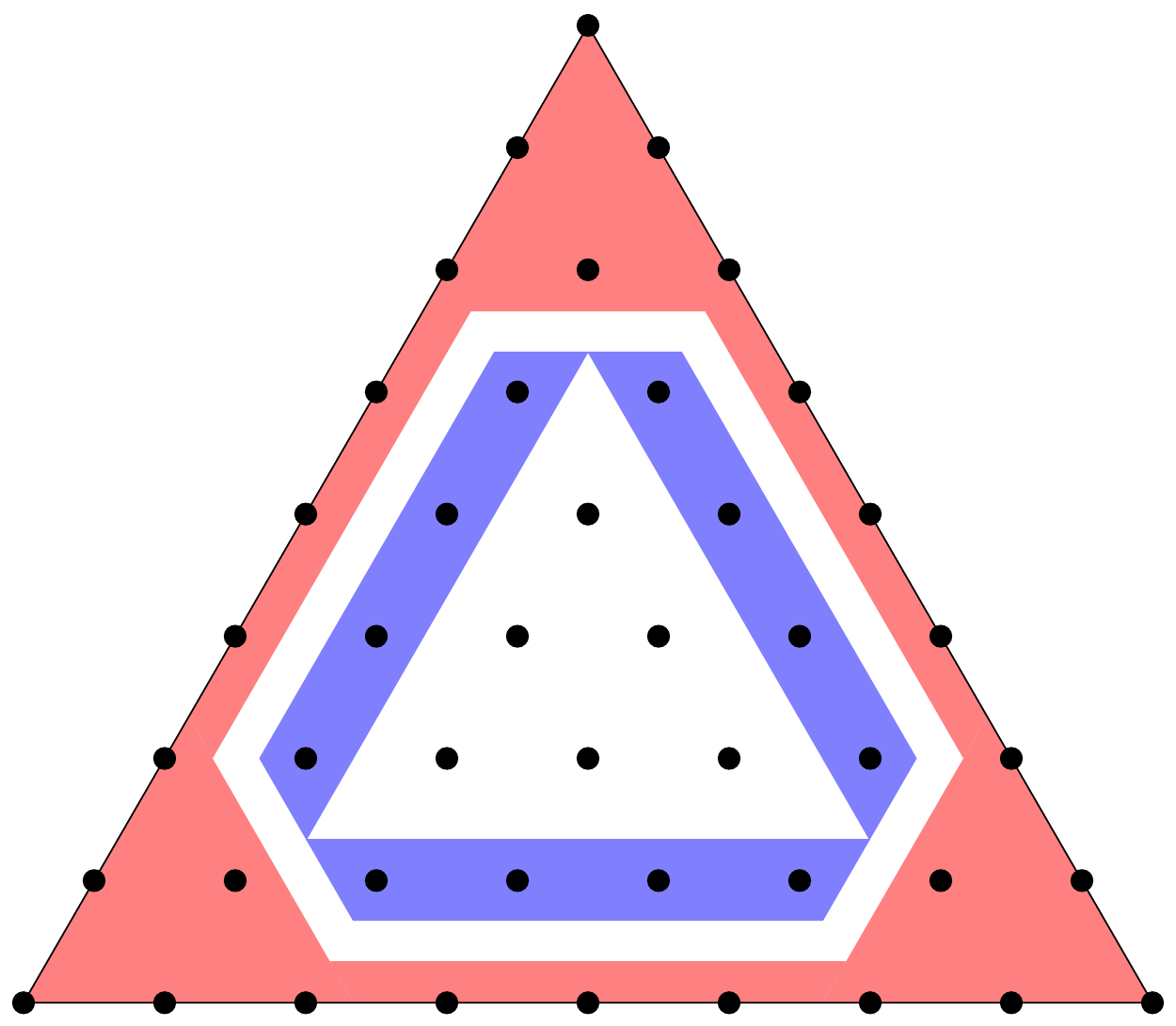}
\end{overpic}
}
\caption{Control points determined by the interpolation conditions for two different types of patches are depicted with black dots. Control points determined from procedures in Section 5.1 and Section 5.2 lie in red area, the ones determined in Section 5.3 lie in blue area.}
\label{fig:S5_S8_triangle}
\end{figure}

In this section we focus on the construction of control points of the sought interpolant, separated into three subproblems:
\begin{itemize}
\item interpolation of tangent planes at vertices (Section 5.1),
\item interpolation of normal curvature forms at vertices  (Section 5.2),
\item approximation of tangent planes at edge midpoints  (Section 5.3).
\end{itemize}
Control points influenced by the interpolation conditions are depicted in Fig.~\ref{fig:S5_S8_triangle}. Algorithms in Sections 5.1 and 5.2 determine boundary control points of the patches and control points near the triangle vertices (control points in red area in Fig.~\ref{fig:S5_S8_triangle}). In Section 5.3 it is explained how to set the remaining control points that influence $C^1$/$G^1$ contacts between patches (control points in blue area in Fig.~\ref{fig:S5_S8_triangle}).

Since the geometric interpolation conditions \corr{would} not set the control points $\bfm c_{\bfm i}$ uniquely, we use the remaining degrees of freedom to obtain well distributed control points by employing control points $\bfm c_{\bfm i}^{\bowtie}$ 
of the referential spline approximant.

In all three cases control points $\bfm c_{\bfm i}^{\bowtie}$ will be projected onto tangent planes.
To achieve $C^1$ smoothness at the vertices, a correction of points will be computed (Section~5.1). The correction will only be needed in the case of $C^1$ quintic patches since in $G^1$ octic case \corr{the smoothness conditions are not directly connected to the} underlying triangulation. To achieve $C^2$ smoothness conditions at the vertices, a similar correction of points will be applied (Section~5.2). In this case the correction will also be \corr{enforced} for the $G^1$ approximant (the projected points in Section~5.1 define a local triangulation that needs to be put into consideration when dealing with $C^2$ conditions).

\subsection{Tangent plane interpolation and minimizing \corrPD{ring of} $\mathcal R_1(v)$}\label{subs:TangPlane}

At every patch vertex $v$ we would like to interpolate a prescribed point $\bfm P$ and the associated tangent plane $\Pi$, defined by the point $\bfm P$ and a normal vector $\bfm n$. To satisfy the first condition we set $\mathcal D_0(v) = \{\bfm P\}$. To interpolate the plane $\Pi$, the constraints
\begin{align}\label{eqn:tangCond}
\langle \bfm c_{\bfm i}-\bfm P, \bfm n \rangle=0, \qquad \bfm c_{\bfm i}\in \mathcal R_1(v),
\end{align}
must be satisfied. The points in $\mathcal D_1(v)$ are connected by smoothness conditions (see Thm.~\ref{thm:smoothPoint}). Hence if we assign positions of the two control points in $\mathcal R_1(v)$ for one patch, then the \corr{remaining ones} in $\mathcal R_1(v)$ are uniquely determined by the $C^1$ continuity conditions. Therefore, the above restrictions give a 4-parametric family of control points.

Since the interpolation conditions are not sufficient to uniquely determine the set $\mathcal D_1(v)$, the remaining degrees of freedom will be used so that the points $\mathcal D_1(v)$ will be close to projected points, obtained from the referential points $\mathcal D_1^{\bowtie}(v)$. 

\begin{figure}[!htb]
\centering
\begin{overpic}[trim=0.5cm 3.5cm .5cm 2cm, clip=true, width=9cm]{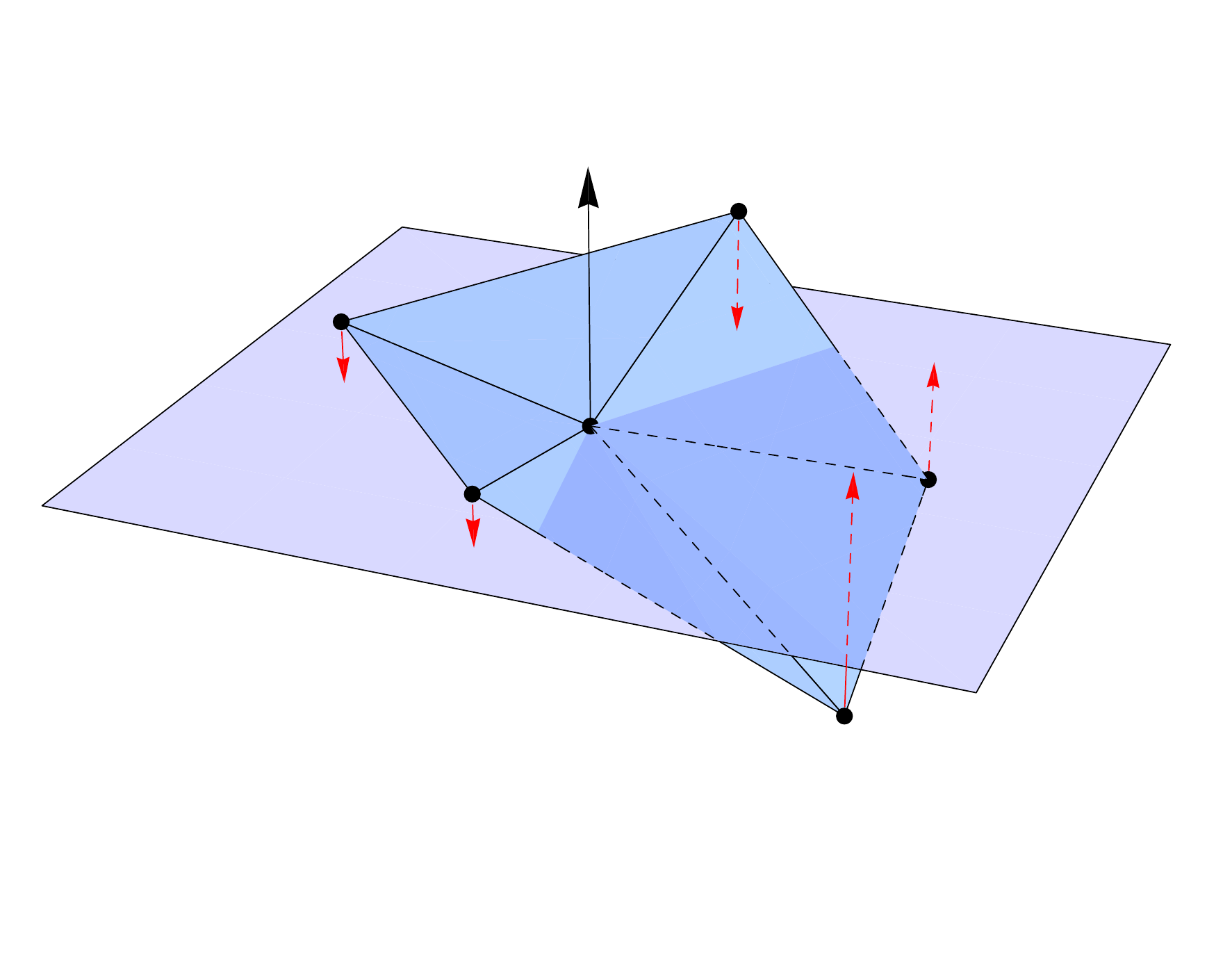}
\put(73,1){\small $\bfm c_1^{\bowtie}$}
\put(80,22){\small $\bfm c_2^{\bowtie}$}
\put(64,46){\small $\bfm c_3^{\bowtie}$}
\put(21,33){\small $\bfm c_4^{\bowtie}$}
\put(33,17){\small $\bfm c_5^{\bowtie}$}
\put(42,25){\small $\bfm P$}
\put(44,36){\small $\bfm n$}

\usetikzlibrary{positioning}
\put(8,23){\small \tikz[baseline] \node at (0,0) [xslant=1.25, yslant=-.15, yscale=.52, xscale=1] {$\Pi$};}
\end{overpic}
\caption{Control points in $\mathcal R_1^{\bowtie}(v)$. Red arrows indicate orthogonal projections of the points onto the plane $\Pi$.}
\label{fig:ordering}
\end{figure}
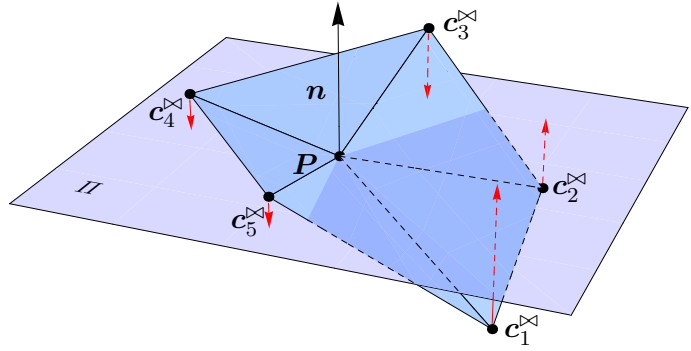

Presume that points in the sets $\mathcal R_1(v)$ and $\mathcal R_1^{\bowtie}(v)$ are denoted by $\bfm c_\ell$ and $\bfm c_\ell^{\bowtie}$ for $\ell = 1,2,\dots,n$, respectively. Let the elements in both sets have the same ordering that corresponds to ordering of vertices around the central vertex in \corrP{the domain} cell (see Fig.~\ref{fig:triangs}(b) and Fig.~\ref{fig:ordering}).
Furthermore, let us presume geometric restrictions
\begin{align}\label{eqn:notParallel}
(\bfm c_\ell^{\bowtie} - \bfm P) \nparallel \bfm n, \qquad \ell=1,2,\dots,n.
\end{align}
The points $\{\bfm c_{\ell}^{\bowtie}\}$ are projected onto $\Pi$ in the direction $\bfm n$ -- they will be denoted by $\big\{\bfm c_\ell^{\rm( p)} \big\}$. 

\begin{figure}[!htb]
\centering
\begin{overpic}[width=7cm]{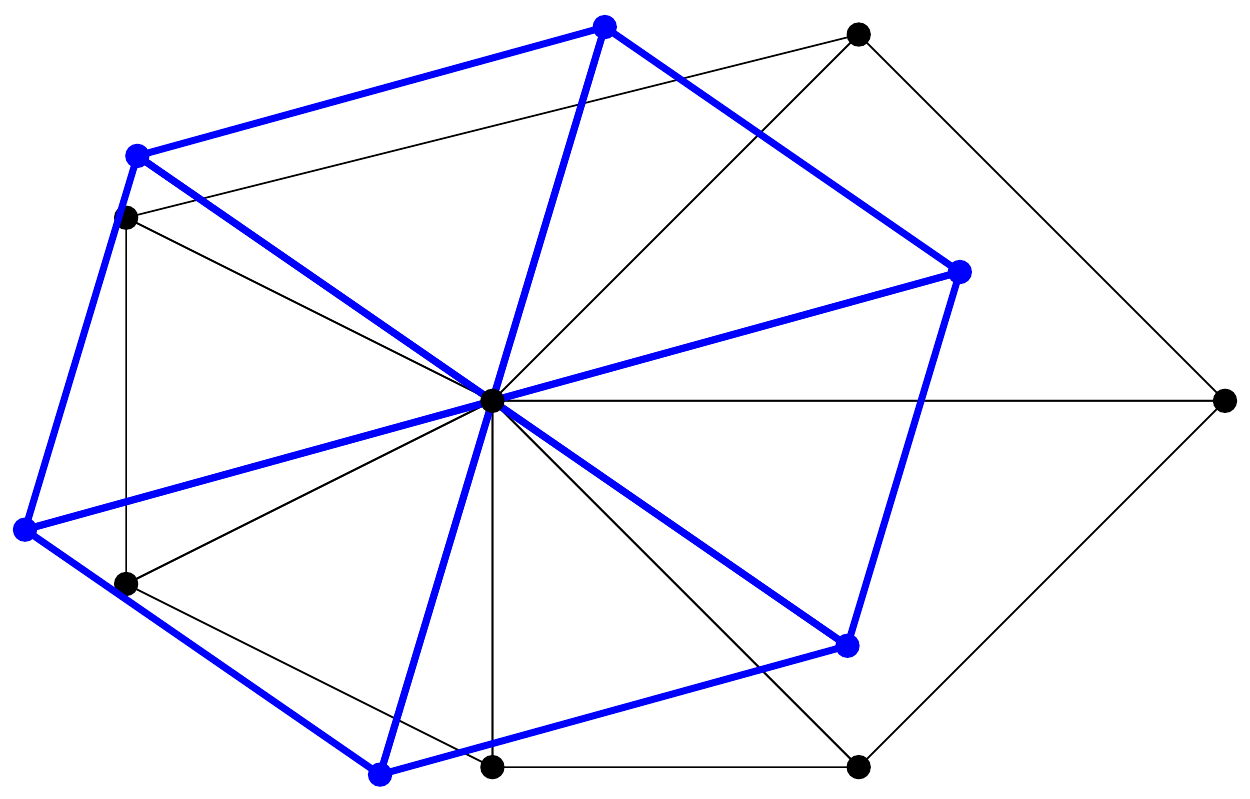}

\put(45,28.5){\small $\bfm P$}

\put(80,41.5){\small $\bfm c_1$}
\put(52,61){\small $\bfm c_2$}
\put(71,11){\small $\bfm c_n$}

\put(101,30.5){\small $\bfm c_1^{(\rm p)}$}
\put(72,60){\small $\bfm c_2^{(\rm p)}$}
\put(73,1){\small $\bfm c_n^{(\rm p)}$}

\put(1,1){\small $\Pi$}

\end{overpic}
\caption{Points $\bfm c_\ell$ of the minimizing \corrPD{ring} are determined in such a way that they are as close to points $\bfm c_\ell^{(\rm p)}$ as possible, i.e., they minimize the functional $\varphi$ in \eqref{eqn:phiFunc}. 
}
\label{fig:minCell}
\end{figure}
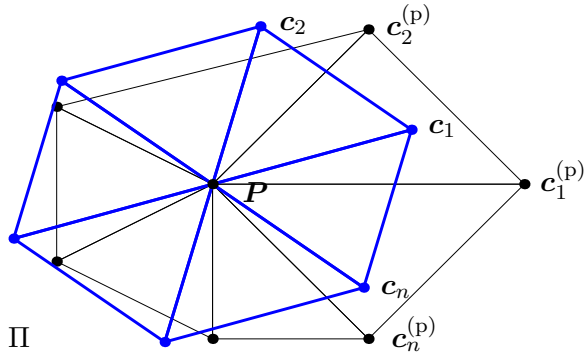

If we would set $\bfm c_\ell := \bfm c_\ell^{(\rm p)}$, $\ell=1,2,\dots,n$, the spline $\bfm s$ would interpolate the plane $\Pi$ but would not be $C^1$ smooth in the neighbourhood of the point $\bfm P$. Therefore, let us find an admissible set of control points $\bfm c_\ell$ with respect to the smoothness conditions that is relatively close to the projected points $\bfm c_\ell^{(\rm p)}$. We would like to solve the least squares minimization problem
\begin{align}\label{eqn:minCell}
\min_{\{\bfm c_1,\,\bfm c_2\}} \varphi ((\bfm c_\ell)_{\ell=1}^n ),
\end{align}
where the functional $\varphi$ measures relative distances between the two sets of points,
\begin{align}\label{eqn:phiFunc}
\varphi
((\bfm c_\ell)_{\ell=1}^n) := \sum_{\ell=1}^n \frac{\left\|\bfm c_\ell - \bfm c_\ell^{(\rm p)} \right\|^2}{\left\|\bfm c_\ell^{(\rm p)} - \bfm P \right\|^2}.
\end{align}
Note that by Thm.~\ref{thm:smoothPoint} the control points $\bfm c_{\ell}$ are connected by $C^1$ smoothness conditions at the vertex,
\begin{align}\label{eqn:mincellCond}
\bfm c_\ell = \left\langle v_\ell(\tau_{\ell-2}), (\bfm P , \bfm c_{\ell-2}, \bfm c_{\ell-1}) \right \rangle, \qquad \ell\ge 3.
\end{align}
Here $v_\ell$ is the vertex that corresponds to the point $\bfm c_\ell$ and $\tau_\ell$ the triangle that corresponds to the points $\bfm c_\ell,\,\bfm c_{\ell+1}$. 
The problem \eqref{eqn:minCell} can be written as a normal equation and it has a unique solution. We call the optimal set of points $\mathcal R_1(v)$ a \emph{minimizing \corrPD{ring}}. An example is shown in Fig.~\ref{fig:minCell}. 

%
%
%

The control points of the minimizing \corrPD{ring} satisfy the tangent plane conditions \eqref{eqn:tangCond}.
\begin{prop}
Let $\mathcal M:=\{\bfm c_\ell \}_{\ell = 1}^n$ be the minimizing \corrPD{ring} of $\mathcal R_1(v)$. Then the points of $\mathcal M$ lie on the plane $\Pi$.
\end{prop}



If no underlying domain triangulation is given when constructing the $G^1$ smooth approximant the procedure to correct the positions of control points $\bfm c_\ell^{\rm (p)}$ by computing the minimizing \corrPD{ring} is omitted. The points $\bfm c_\ell^{\rm (p)}$ themselves define a local domain triangulation which will be \corr{used} in Section 5.2.


\subsection{Normal curvature interpolation and minimizing \corrPD{ring of} $\mathcal R_2(v)$}\label{subs:curvInterp}

In this subsection we presume that the set of points $\mathcal D_1(v)$, $v\in\mathcal V$, is already fixed (e.g., it is determined by the procedure in Section~\ref{subs:TangPlane}) so that the spline $\bfm s$ interpolates a point $\bfm P$ and a tangent plane, defined by a normal $\bfm n$, at $v$. The remaining points in $\mathcal D_2(v)$ will be used to interpolate a given normal curvature form at $v$.  The form is described by a set
\begin{align}\label{eqn:curvQuad}
\{\bfm u_{1}^*, \bfm u_{2}^*, \kappa_{1}, \kappa_{2}\},
\end{align}
where $\bfm u_{\ell}^*$ and $\kappa_{\ell}$ are the principal directions and the corresponding normal curvatures of a surface at $v$.
A well known property from differential geometry states that the normal curvature $\kappa_{\textrm n}(\bfm u)$ of the spline $\bfm s$ in direction $\bfm u$, $\|\bfm u\| =1$, is
\begin{align*}
\kappa_{\textrm n}(\bfm u)= \kappa_{1} \langle \bfm u, \bfm u_{1}^* \rangle^2 + \kappa_{2} \langle \bfm u, \bfm u_{2}^* \rangle^2.
\end{align*}
The presumption $\bfm s \in C^2(v)$ ensures the consistency of the curvature form of the neighbouring patches.
Before dealing with the construction of control points, we need the following lemma that states the connection between the normal curvatures and the control points.

The points in $\mathcal R_2(v)$ are connected by $C^2$ smoothness conditions (see Thm.~\ref{thm:smoothPoint}). Note that $\mathcal D_1(v)$ defines a local domain triangulation needed for $C^2$ conditions if no underlying domain triangulation is given. If we fix the three control points of one of the surrounding patches in $\mathcal R_2(v)$, then the rest in $\mathcal R_2(v)$ are uniquely determined by $C^2$ continuity constraints. The above conditions define a 6-parametric family of control points (3 out of 9 degrees of freedom are determined from the normal curvature form).

\begin{figure}[!htb]
\centering
\begin{overpic}[width=7cm]{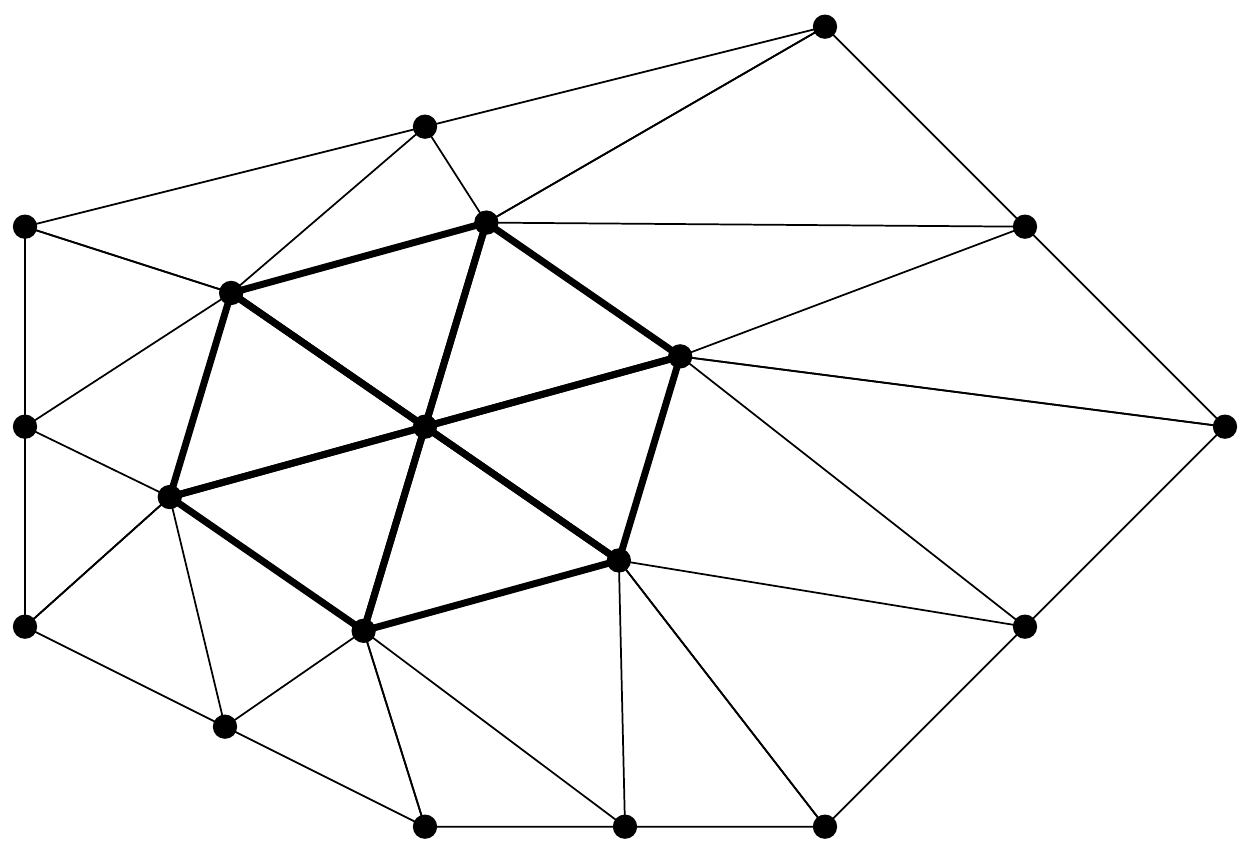}

\put(30,36.5){\small $\bfm P$}

\put(58,39){\small $\bfm c_1$} \put(63,39){\scriptsize $= \bfm c_{\eta(n+1)}$}
\put(43,51.3){\small $\bfm c_2$} \put(49,51.3){\scriptsize $= \bfm c_{\eta(n+3)}$}
\put(52,21.5){\small $\bfm c_n$}

\put(100.5,32.5){\small $\bfm c_{n+1}$}
\put(84,48.5){\small $\bfm c_{n+2}$}
\put(68.5,64.5){\small $\bfm c_{n+3}$}

\put(83.5,17){\small $\bfm c_{n+n'}$}

\end{overpic}
\caption{Points in $\mathcal R_1(v)$ and $\mathcal R_2(v)$ are indexed with the same ordering.}
\label{fig:minCell2}
\end{figure}

The remaining 6 parameters will be obtained from the minimizing \corrPD{ring}. Let us presume that the control points $\mathcal R_1(v)=:\{\bfm c_\ell \}_{\ell=1}^n$ of $\bfm s$ are ordered as in Section~\ref{subs:TangPlane}. Let points in $\mathcal R_2(v)=:\{\bfm c_\ell\}_{\ell=n+1}^{n+n'}$ and $\mathcal R_2^{\bowtie}(v)=:\{\bfm c^{\bowtie}_\ell\}_{\ell=n+1}^{n+n'}$ be indexed with the same ordering as $\mathcal R_1(v)$ (see Fig.~\ref{fig:minCell2}). Here, $n' = 2n$ if $v$ is interior and $n' = 2n-1$, otherwise.


The points $\bfm c_\ell^{\bowtie}$ are projected in the direction of $\bfm n$,
\begin{align}\label{eqn:projPoint2}
&\bfm c_\ell^{(\rm p)} := \bfm c_\ell^{\bowtie} - (\langle \bfm c_\ell^{\bowtie} - \bfm P, \bfm n \rangle - k_\ell)\, \bfm n,
\qquad \ell=n+1,n+2,\dots,n+n',
\end{align}
where
\begin{align}
k_\ell := \left\{
\begin{array}{ll}\label{eqn:curvTerm}
\displaystyle \frac{d}{d-1} \left \| \bfm c_{\eta(\ell)} -  \bfm P \right \|^2 \kappa_{\textrm n}(\bfm u_\ell), & \textrm{if } \ell-n \textrm{ is odd},\\
\displaystyle \frac{2d}{d-1} \left \| \corr{\frac{1}{2}( \bfm c_{\eta(\ell-1)} + \bfm c_{\eta(\ell+1)})} -  \bfm P \right \|^2 \kappa_{\textrm n}(\bfm u_\ell) \\
\; \qquad - \displaystyle \frac{1}{2} \left \langle \bfm c_{\ell-1}^{(\rm p)} - \bfm c_{\eta(\ell-1)} + \bfm c_{\ell+1}^{(\rm p)}- \bfm c_{\eta(\ell+1)} , \bfm n \right \rangle, & \textrm{if } \ell-n \textrm{ is even},
\end{array}
\right.
\end{align}
and
\begin{align*}
\corr{\bfm c_{\eta(\ell)} := \bfm c_{(\ell - n+1)/2}}, \qquad
\bfm u_\ell &:= \frac{1}{\|\bfm c_{\eta(\ell)} -  \bfm P \|}(\bfm c_{\eta(\ell)} -  \bfm P).
\end{align*}
In \eqref{eqn:projPoint2} we first need to compute the points $\bfm c_\ell^{(\rm p)}$ where $\ell-n$ is odd. 

Setting $\bfm c_\ell := \bfm c_\ell^{(\rm p)}$, $\ell = n+1,n+2,\dots,n+n'$, would result in a spline $\bfm s$ that interpolates the normal curvature form at $v$ but is not $C^2$ smooth at $v$. Therefore, we need to find a set of points that satisfies the smoothness constraints and is close to the points $\bfm c_\ell^{(\rm p)}$. Hence, we use the functional $\varphi$, introduced in \eqref{eqn:phiFunc}, and  solve the minimization problem
\begin{align}\label{eqn:minCell2}
\min_{\{ \bfm c_{n+1},\, \bfm c_{n+2},\, \bfm c_{n+3} \}} \varphi \left((\bfm c_\ell)_{\ell=n+1}^{n+n'} \right).
\end{align}
The control points $\bfm c_\ell$, $\ell > n+3$, are uniquely set from \corr{$\bfm c_{n+1}, \bfm c_{n+2}, \bfm c_{n+3}$} by the corresponding $C^2$ smoothness conditions at $v$ (see Thm.~\ref{thm:smoothPoint}).

As in the tangent plane interpolation problem, we are left to verify that the control points in the minimizing \corrP{ring} satisfy the normal curvature interpolation conditions. Let $\Pi_\ell$ denote a plane defined by a point $\bfm c_\ell^{(\rm p)}$ and the normal vector $\bfm n$. Then the curvature constraints are transformed to
\begin{align}\label{eqn:curvatureCond}
\bfm c_\ell \in \Pi_\ell, \qquad \ell = n+1,n+2,\dots, n+n'.
\end{align}

\begin{lem} \label{lem:overlapCond} Let points in the set $\mathcal D_2(v)$ satisfy $C^2$ smoothness conditions at $v$. \corr{If there exists $k \in 2(\NN \cup \{0\})$ such that $\bfm c_\ell \in \Pi_\ell$ for $\ell=n+k+1,n+k+2, n+k+3$ (i.e., a triple of points $\bfm c_\ell$ that correspond to the same patch), then \eqref{eqn:curvatureCond} holds true.}
\end{lem}

%

\begin{prop}
Let $\mathcal M:=\{\bfm c_\ell\}_{\ell=n+1}^{n+n'}$ be the minimizing \corrPD{ring} for \eqref{eqn:minCell2}. Then the points of $\mathcal M$ satisfy conditions \eqref{eqn:curvatureCond}.
\end{prop}

\subsection{Tangent plane approximation at midpoint of an edge}\label{subs:TangPlane2}

In the last part of the section we will analyze the problem on how to determine the remaining control point $\bfm v_2$ of the transversal vector function $\bfm v$ in order to approximate a given tangent plane.

Let us presume that steps in Sections~5.1 and 5.2 are already applied and that control points $\bfm v_\ell$, $\ell = 0,1,3,4$, are appropriately chosen \corr{(see 
Fig.~\ref{fig:g1_pog})}.

%

Let $\bfm n_\Pi$ denote the normal vector of the plane $\Pi$ that we would like to approximate at the edge midpoint $\bfm p^{[\tau_2]}(v_0/2+v_2/2)$. Since the tangent vector $\bfm b(1/2)$ already fixes one direction of the tangent plane of patches at the boundary, we can only approximate $\Pi$.
To obtain the best approximating tangent plane of $\Pi$ (denoted by $\Pi^\star$), the tangent plane normal $\bfm n$ of the patch should be set in such a way that $\| \bfm n - \bfm n_\Pi \|$ is minimal. 
Thus, $\bfm n$ is set as orthogonal projection of $\bfm n_\Pi$ onto plane defined by the point $\bfm p^{[\tau_2]}(v_2/2+v_0/2)$ and the normal in the direction of $\bfm b(1/2)$.

Let us decompose $\bfm v$ and its control points into two parts $\bfm v=: \bfm v^{(\bfm n)} + \bfm v^{(\Pi^\star)}$, the first part is a component in the direction of the normal $\bfm n$, the other components is orthogonal to $\bfm n$. The interpolation of the tangent plane $\Pi^\star$ thus reads $\left\langle \bfm v(1/2), \bfm n  \right\rangle = 0$, hence
\begin{align}\label{eqn:v2normal}
\bfm v_2^{(\bfm n)} = \left\langle  \bfm v_2, \bfm n  \right\rangle = 
-\left\langle \sum_{\ell=0,1,3,4} \bfm v_\ell\, \frac{B_\ell^4(1/2)}{B_2^4(1/2)}, \bfm n  \right\rangle =
-\sum_{\ell=0,1,3,4} \bfm v_\ell^{(\bfm n)}\, \frac{B_\ell^4(1/2)}{B_2^4(1/2)}.
\end{align}
Components of $\bfm v$ orthogonal to $\bfm n$ do not influence the interpolation conditions. We set $\bfm v_2^{(\Pi^\star)}$ so that $\bfm v^{(\Pi^\star)}$ is a cubic polynomial:
\begin{align}\label{eqn:v2ort}
\bfm v_2^{(\Pi^\star)} = \frac 1 6 \left( -\bfm v_0^{(\Pi^\star)} + 4\,\bfm v_1^{(\Pi^\star)} + 4\,\bfm v_3^{(\Pi^\star)} -  \bfm v_4^{(\Pi^\star)}\right).
\end{align}

Once \corrP{the} vector function $\bfm v$ is fixed the remaining control points that influence $C^1$/$G^1$ conditions (control points in the blue region in Fig.~\ref{fig:S5_S8_triangle}) are determined via computing vectors \corrP{$\bfm d'_\ell,\bfm e'_\ell$} from \eqref{eqn:G1_cond1} and \eqref{eqn:G1_cond2}.\\

\section{Numerical examples}

Let us conclude the paper by some numerical examples. Our $C^1$ quintic and $G^1$ octic schemes  are tested by approximating a torus and a more general free-form surface. The results are compared \corrP{against three} $G^1$ interpolation schemes. \corrP{The first one is a} scheme by Shirman and S\'equin (SS), \corrP{a quartic 3-splitting scheme, which interpolates points and the corresponding tangent planes at the vertices} \cite{ShirmanSequin, ShirmanSequin_error}. \corrP{The scheme} by Hahmann and Bonneau (HB) \corrP{is a quintic 4-splitting method, which interpolates points at the vertices \cite{Hahmann-Bonneau-G1-triang}. In our tests we} use shape parameters that were also used by the authors and seem to produce the best results: $\alpha=1,\beta=0.1,\gamma_0=-3.7,\gamma_1=4.6,\gamma_2=0.1$. \corrP{The third method is by Tong and Kim (TK) \cite{Tong-Kim-G1-triang}. The patches of degree 7 interpolate points, tangent planes and normal curvatures at the vertices. The remaining degrees of freedom are used to minimize a particular energy functional and the distance to the original surface by applying subsequent data fitting procedures. A basic quantitative comparison between the methods is presented in Tab.~\ref{tab:compareBasic}.}

\begin{table}[!htb]
\centering \small
\begin{tabular}{l|c|c|c|c|c}
\rule{0ex}{2.5ex} & $C^1$ quintic & $G^1$ octic & SS & HB & \corrP{TK}\\
\hline
\rule{0ex}{2.5ex} Polynomial degree & 5 & 8 & 4 & 5 & \corrP{7}\\
\rule{0ex}{2.5ex} \corrP{Micro} patches & \corrP{1} & \corrP{1} & \corrP{3} & \corrP{4} & \corrP{1}\\
\rule{0ex}{2.5ex} Control points & \corrP{21} & \corrP{45} & \corrP{45} & \corrP{84} & \corrP{36}\\
\rule{0ex}{2.5ex} \corrP{Approximation data} & \corrP{27} & \corrP{27} & \corrP{15} & \corrP{9} & \corrP{81}
\end{tabular}
\caption{\corrP{Basic comparison of $G^1$ approximation splines. Number of micro patches, control points and scalar approximation data are counted on a macro patch.}}
\label{tab:compareBasic}
\end{table}

\corrP{In the third numerical example, our scheme is tested against the standard functional Argyris element} by approximating a nonparametric surface. \corr{Short numerical test of the approximation order is done at the end of the paper.}

In all the examples, we fix the transversal vector $\bfm v$ in our interpolation schemes to satisfy conditions \eqref{eqn:simpleVectorV}. The remaining six interior control points $\mathcal C_{\rm i}:=\{\bfm c_{ijk}: i,j,k \geq 2\}$ on every $G^1$ octic spline patch (see unmarked control points in Fig.~\ref{fig:S5_S8_triangle}(b)) are determined by applying a linear combination of all the other control points $\mathcal C_{\rm b}:=\{\bfm c_{\bfm i}: |\bfm i|=8 \}\backslash \mathcal C_{\rm i}$ on the same patch. More precisely, control points $\mathcal  C_{\rm i}$ are defined in such a way that control points $\{\bfm c_{\bfm i}: |\bfm i|=8 \}$ represent a quintic patch if the other control points $\mathcal C_{\rm b}$ are also obtained from that same patch. An alternative that is not explored in this paper would be to use the remaining six points to minimize particular energy functionals (see \cite{DiscrWill}, e.g.) or to interpolate additional points in the interior - such interpolation problems are \corr{unisolvent} \cite{Jaklic-Kanduc-BezierDet, Schumaker-Lai-splines-triangulations-07}.

\subsection{Torus approximation}
\begin{figure}[!htb]
\centering
\begin{overpic}[trim=0cm .4cm 0 .3cm, clip=true, width=5.5cm]{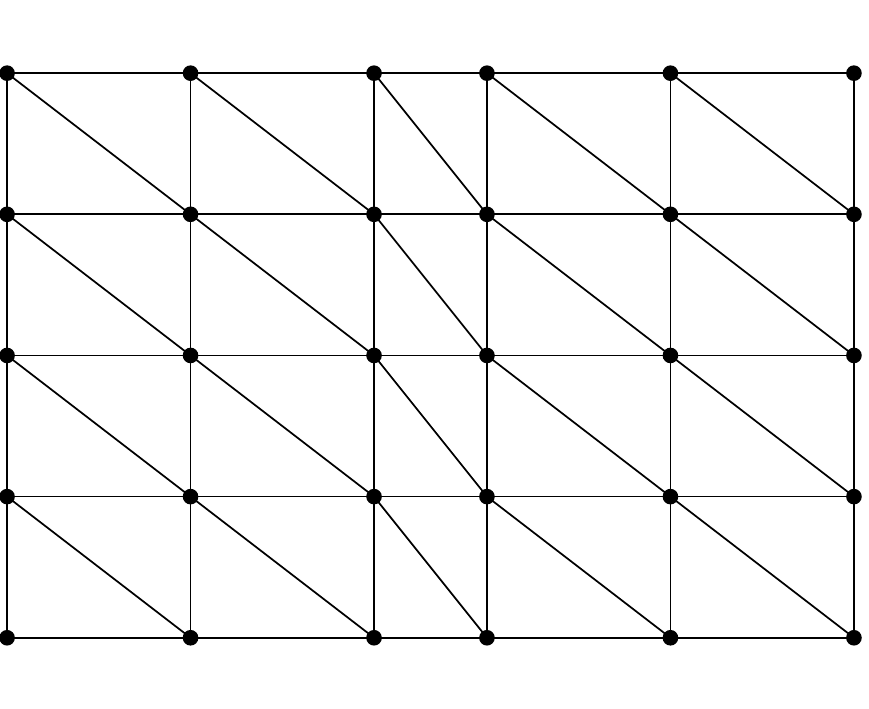}
\put(-2,71){\small $1$}
\put(20,71){\small $2$}
\put(40,71){\small $3$}
\put(57,71){\small $4$}
\put(78,71){\small $5$}
\put(99,71){\small $1$}

\put(-2,-4){\small $1$}
\put(20,-4){\small $2$}
\put(40,-4){\small $3$}
\put(57,-4){\small $4$}
\put(78,-4){\small $5$}
\put(99,-4){\small $1$}

\put(-4,55){\small $6$}
\put(-4,38.5){\small $7$}
\put(-4,22.5){\small $8$}

\put(100,55){\small $6$}
\put(100,38.5){\small $7$}
\put(100,22.5){\small $8$}
\end{overpic}
\caption{Domain triangulation for the torus, where vertices with the same indices and the corresponding edges are identified.}
\label{fig:torusDomain}
\end{figure}
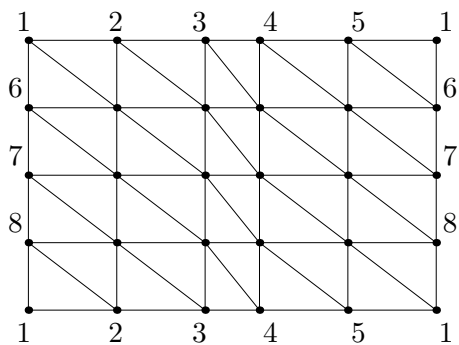

In the first example we approximate a torus with a major radius $R=2$ and 
a minor radius $r=1$.
By identifying boundary vertices and edges of the domain triangulation, we construct a triangulation 
suitable for approximating a torus with $C^1$ smooth splines (Fig.~\ref{fig:torusDomain}).

\begin{figure}[!htb]
\centering
\subfigure[Referential surface]{
\includegraphics[trim=-2cm .0cm -2cm 0cm, clip=true, height=4cm]{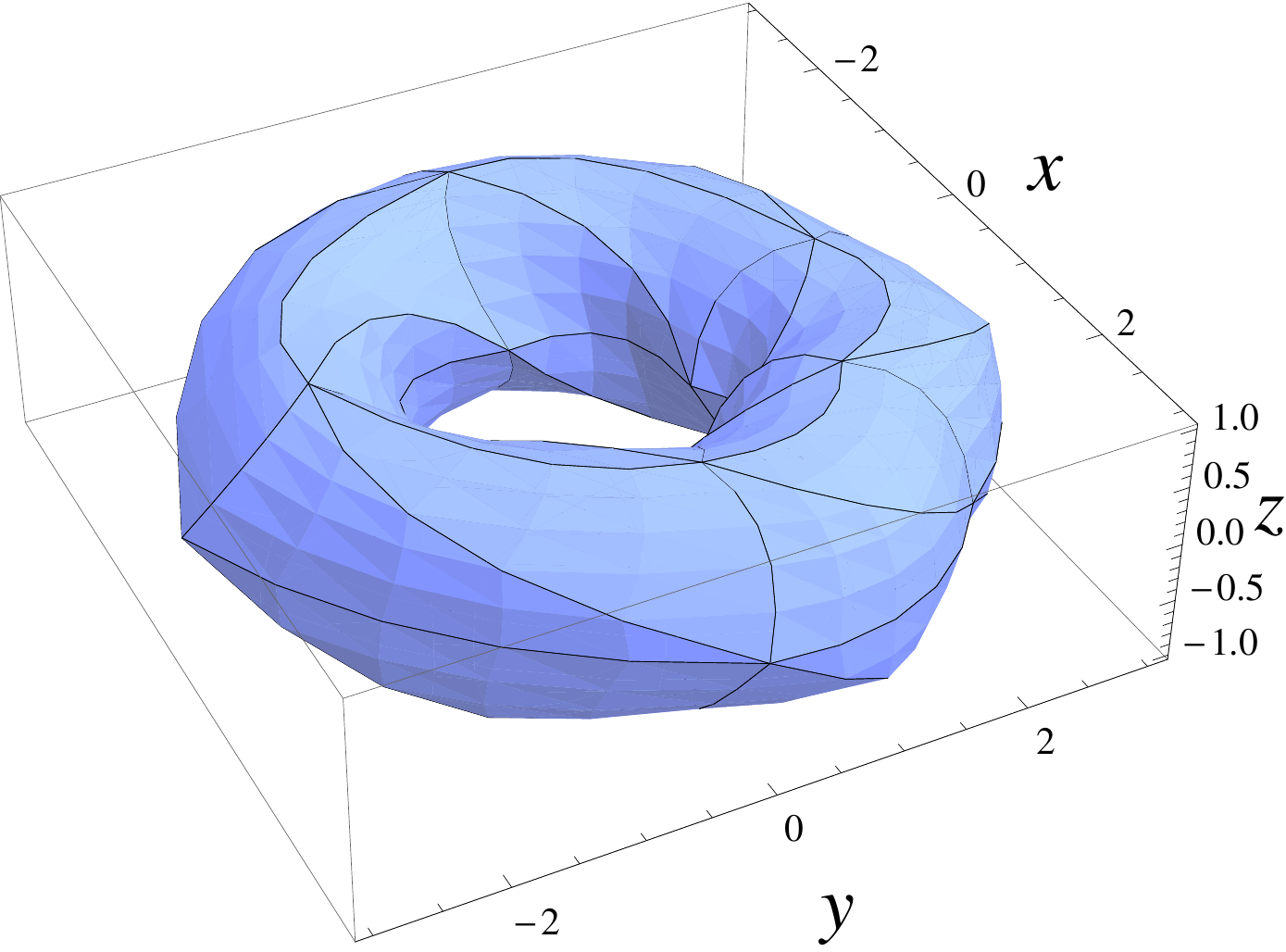}
}
\subfigure[$C^1$ quintic approximant, \corrP{$E_H$:  0.13}]{
\includegraphics[trim=-2cm .0cm -2cm 0cm, clip=true, height=4cm]{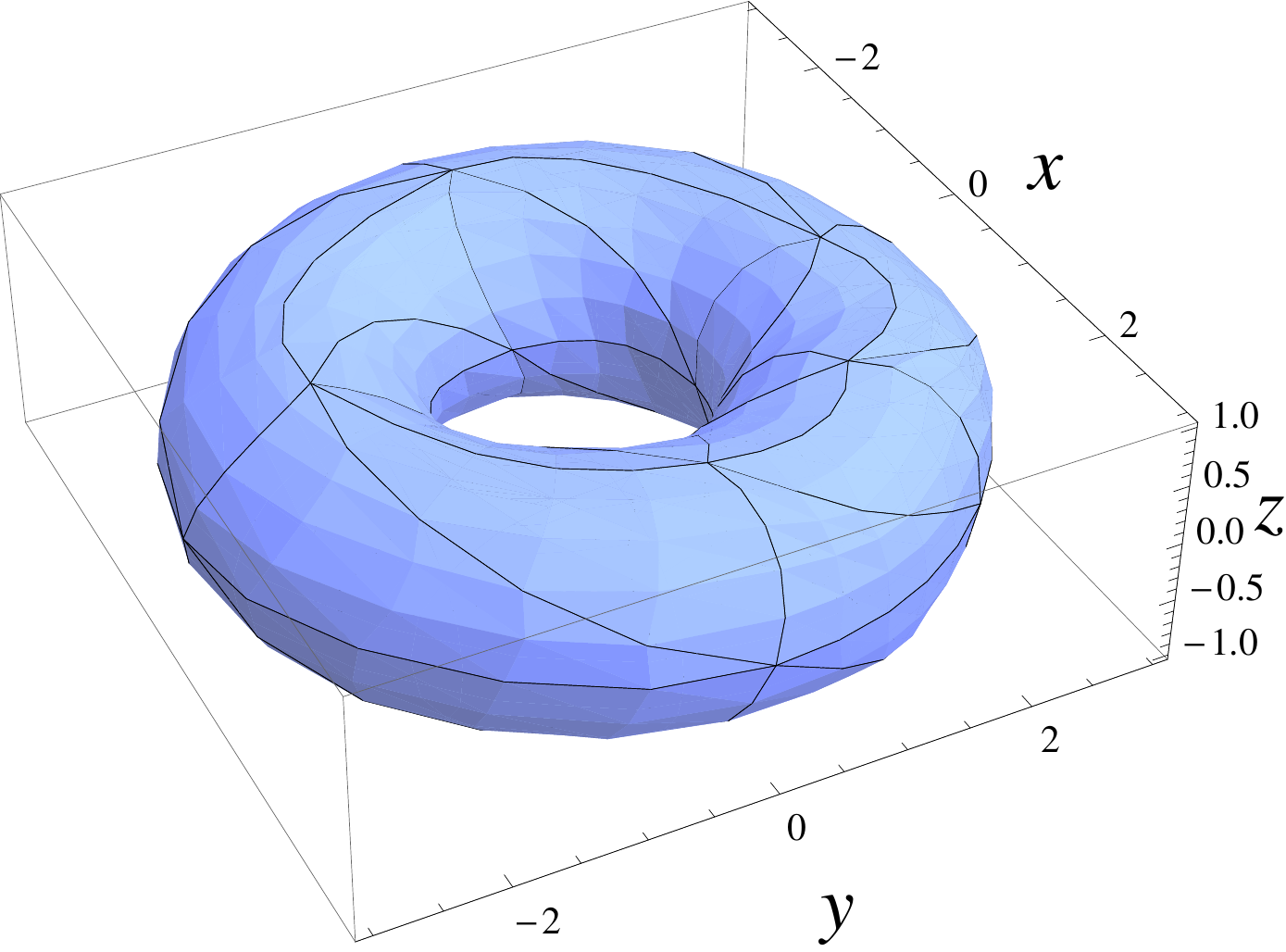}
}
\subfigure[$G^1$ octic approximant, \corrP{$E_H$:  0.13}]{
\includegraphics[trim=-2cm .0cm -2cm 0cm, clip=true, height=4.cm]{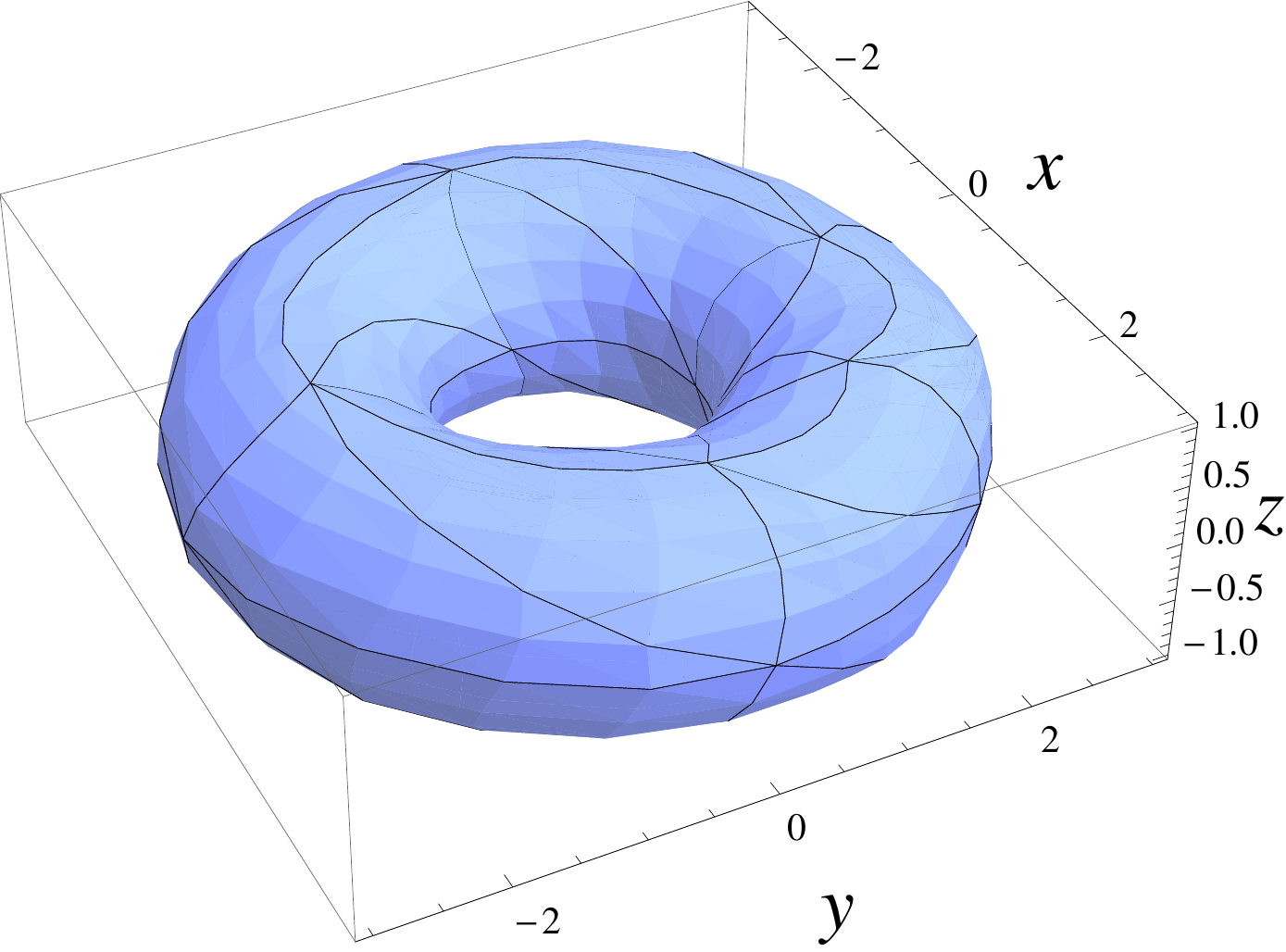}
}
\subfigure[SS approximant, \corrP{$E_H$:  0.45}]{
\includegraphics[trim=-2cm .0cm -2cm 0cm, clip=true, height=4.cm]{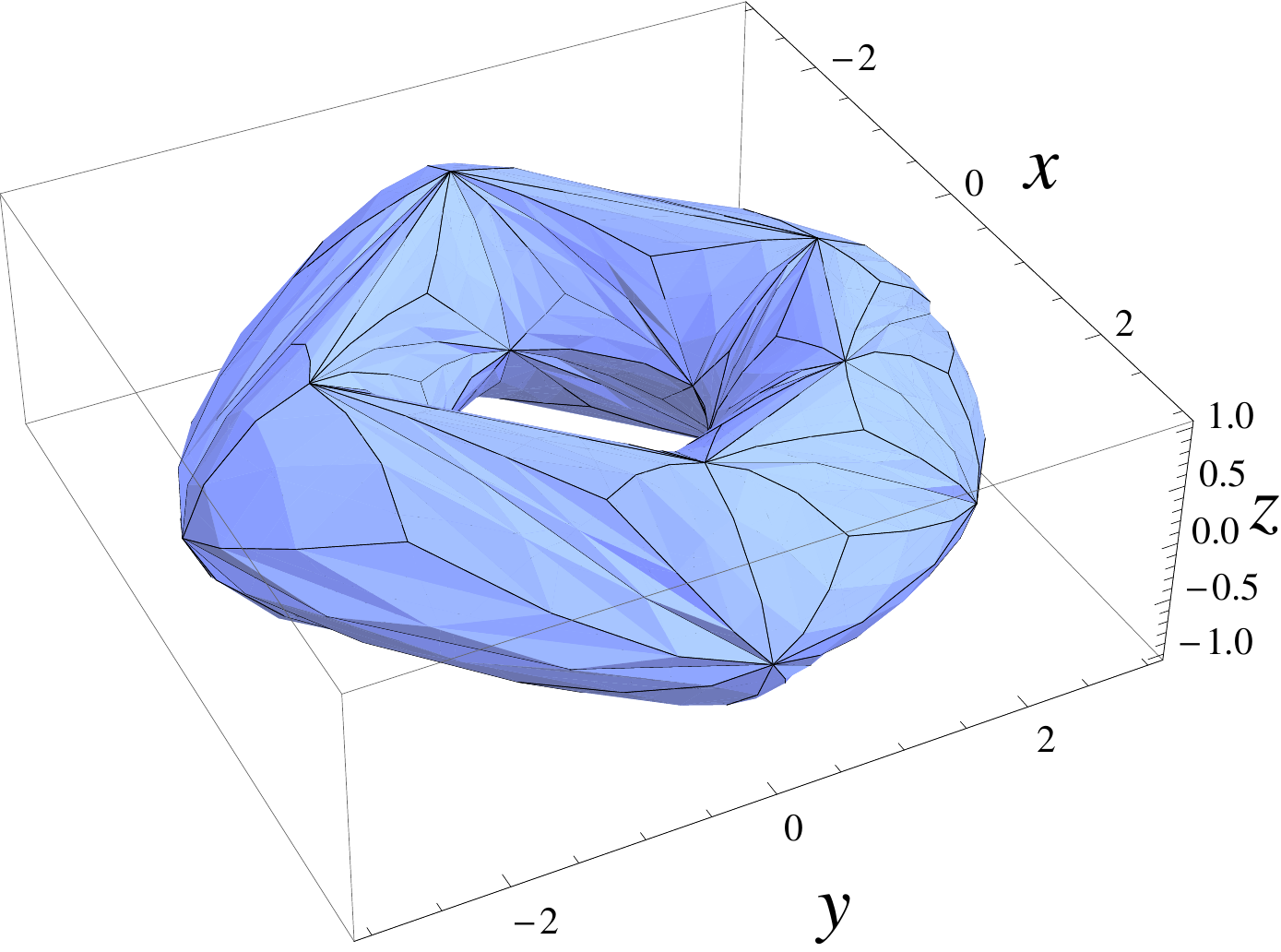}
}
\subfigure[HB approximant, \corrP{$E_H$:  0.59}]{
\includegraphics[trim=-2cm .0cm -2cm 0cm, clip=true, height=4.cm]{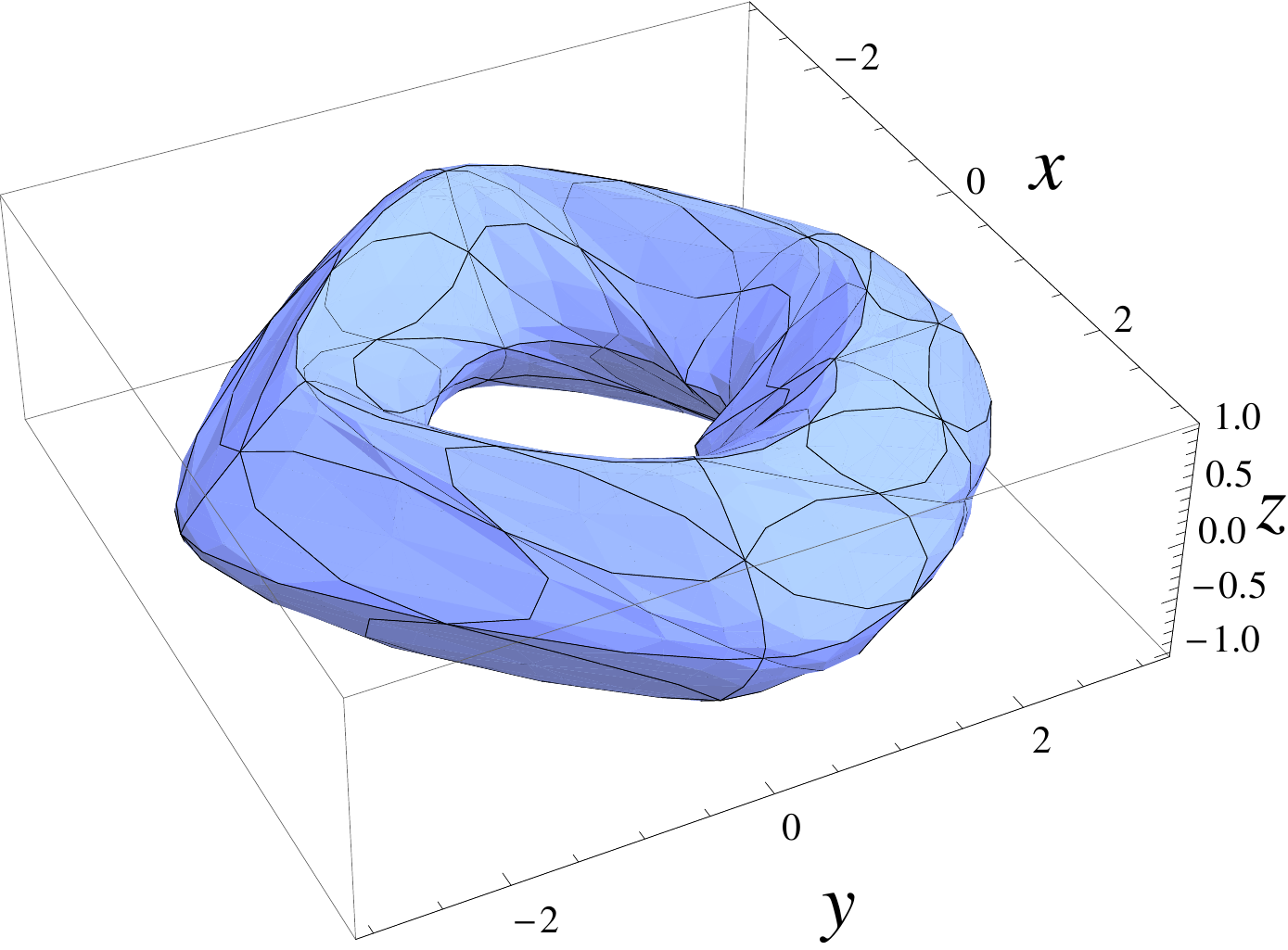}
}
\subfigure[\corrP{TK approximant, $E_H$:  0.26}]{
\includegraphics[height=4.cm]{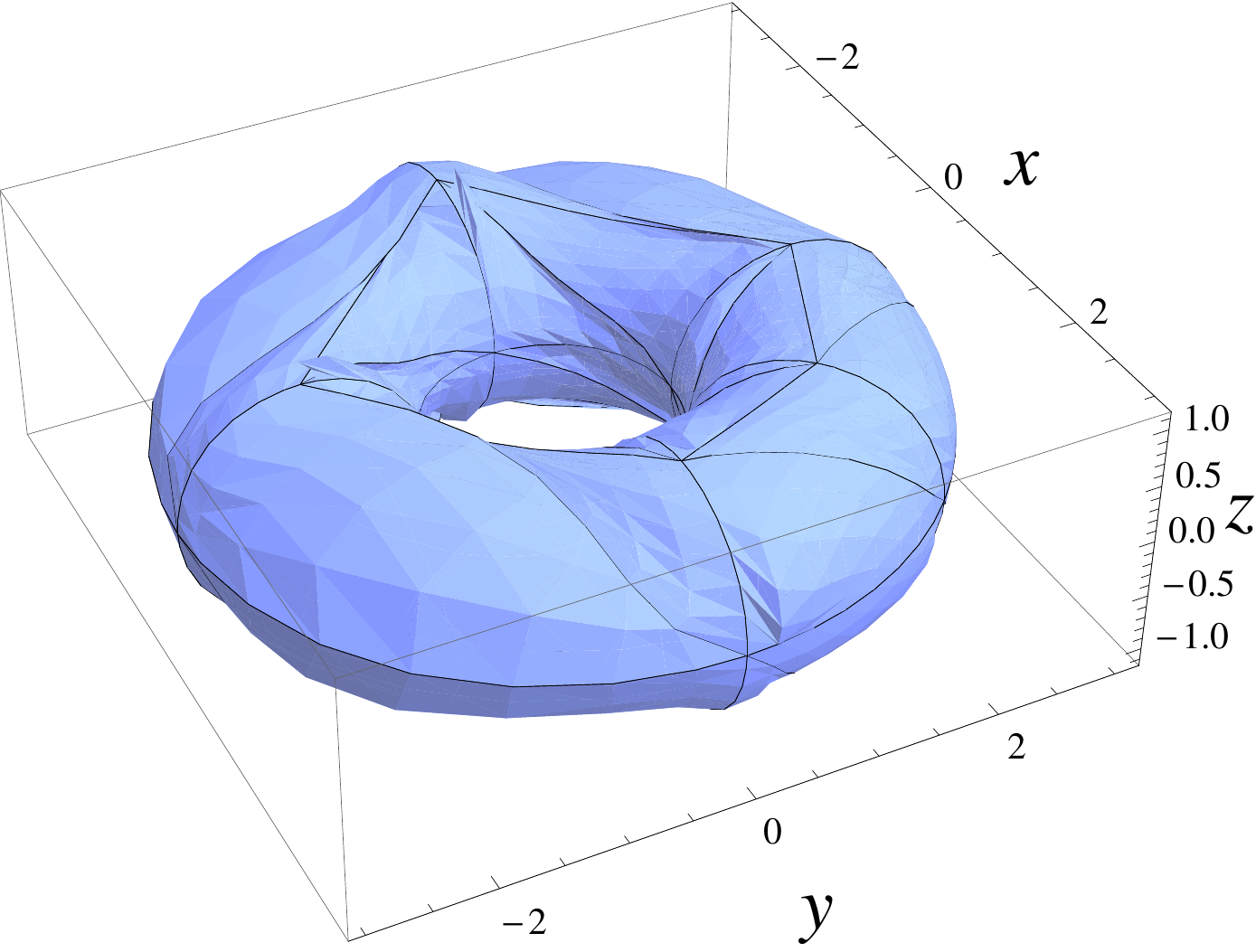}
}
\caption{Approximation of \corrP{the} torus by different interpolation schemes. \corrP{Hausdorff errors ($E_H$) are shown next to the plots.}}
\label{fig:torusApprox}
\end{figure}

The torus is approximated by the introduced $C^1$ quintic and $G^1$ octic splines (see Fig.~\ref{fig:torusApprox}). Referential surface -- an intermediate step to construct the interpolation surface -- is also depicted. To test the quality of the approximants, a comparison is made with \corrP{SS, HB and TK scheme}. 
All the interpolants approximate torus better at the right-hand side segments since the interpolation data are denser in that area. Both of our approximants have smaller Hausdorff errors than the other \corrP{three} schemes. This can be partially justified by \corrP{the fact that our approximants interpolate more data than SS and HB scheme}. Furthermore the surface curvature is apparently better distributed along the spline patches. \corrP{For smaller patches TK scheme produces a surface with small Hausdorff error. Undesired intersections of boundary curves can be observed on the top part of the surface. We have also noticed that the method is sensitive to the input data and to the parameters used in the minimizing processes of the algorithm.}

\subsection{Free-form surface approximation}
\begin{figure}[!htb]
\centering
\subfigure[Surface defined by $\bfm f$]{
\includegraphics[height=4.5cm]{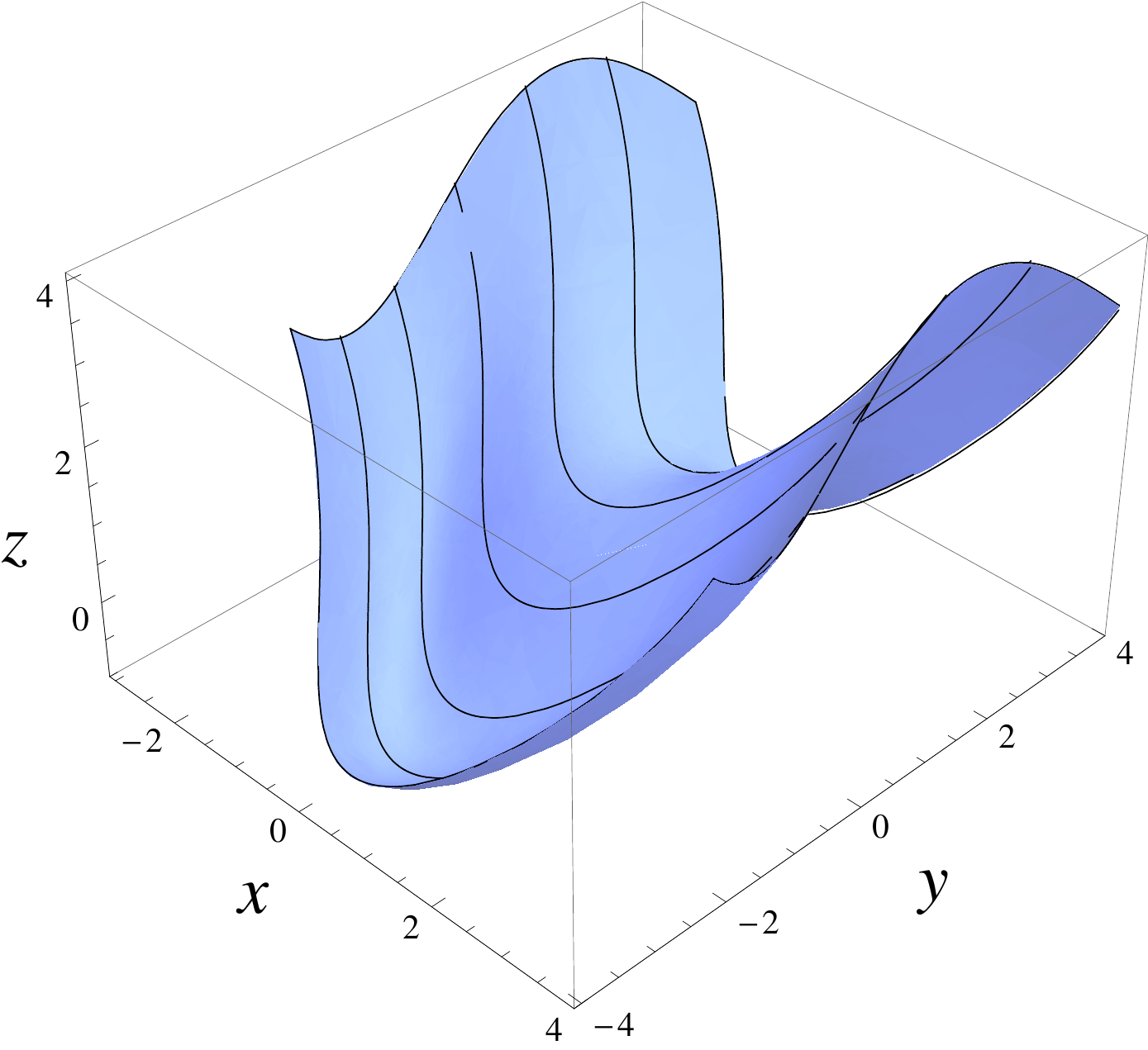}
}
\subfigure[Referential surface]{
\includegraphics[height=4.5cm]{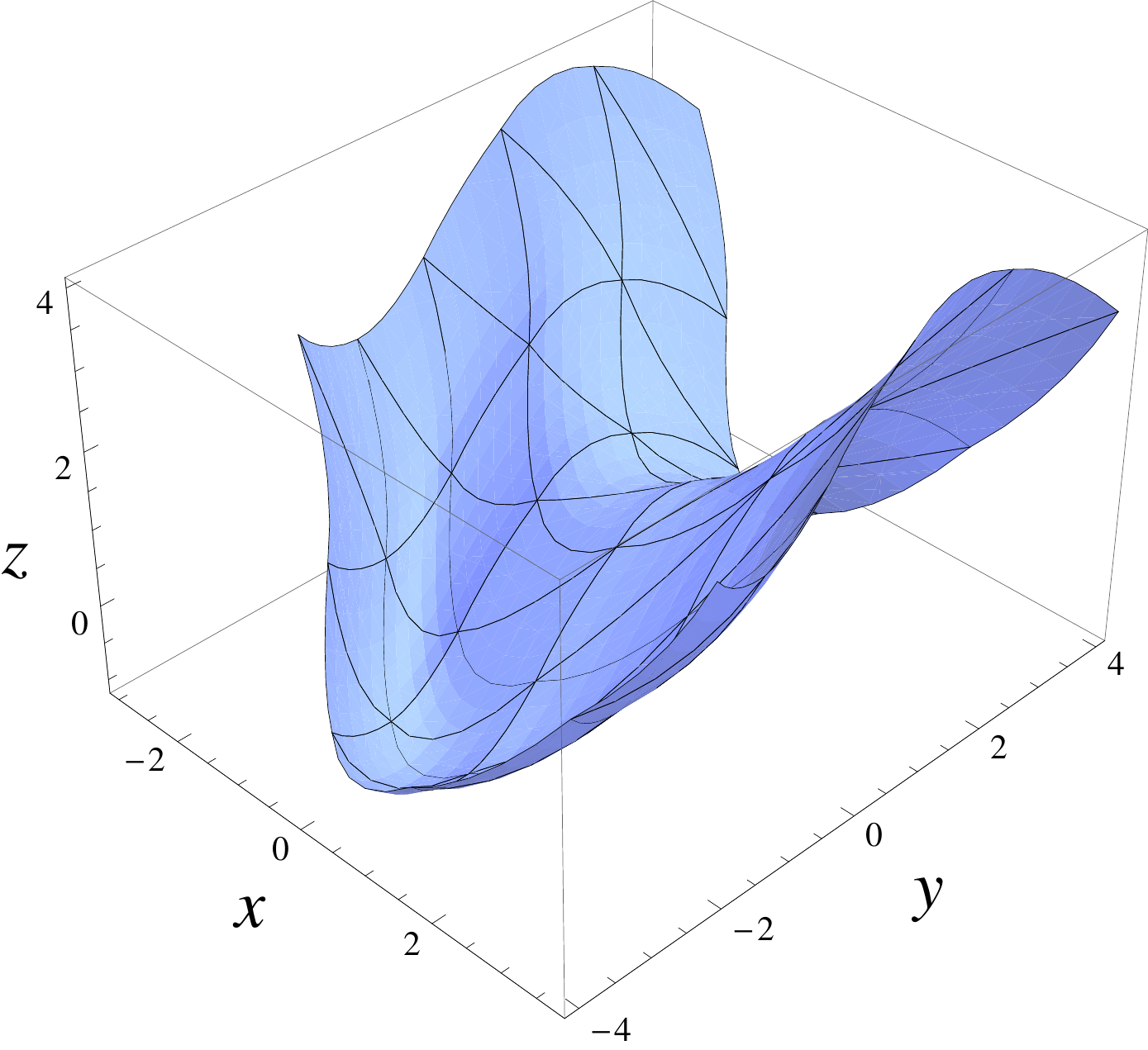}
}
\subfigure[$C^1$ quintic scheme, \corrP{$E_H$:  0.051}]{
\includegraphics[height=4.5cm]{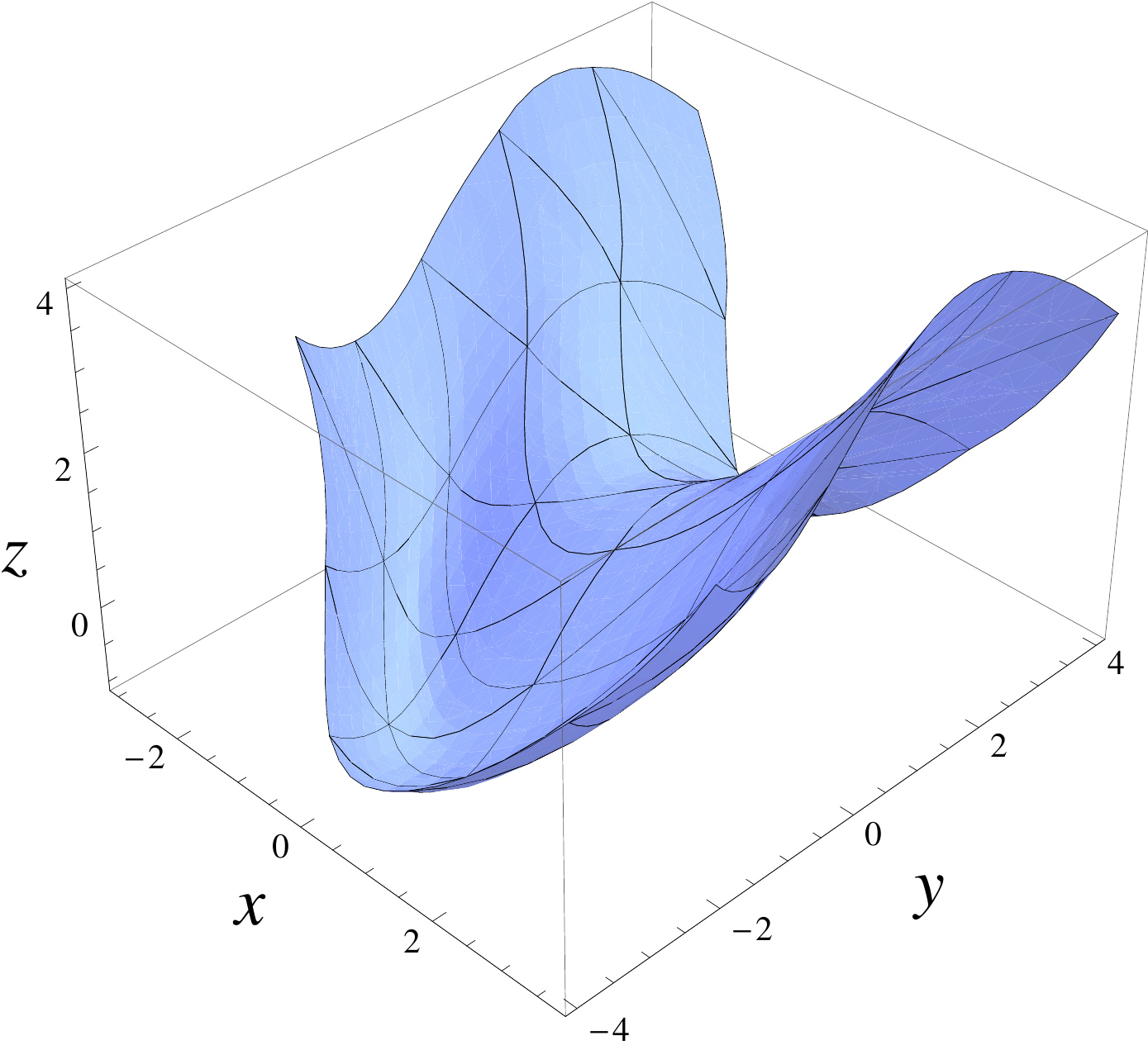}
}
\subfigure[$G^1$ octic scheme, \corrP{$E_H$:  0.050}]{
\includegraphics[height=4.5cm]{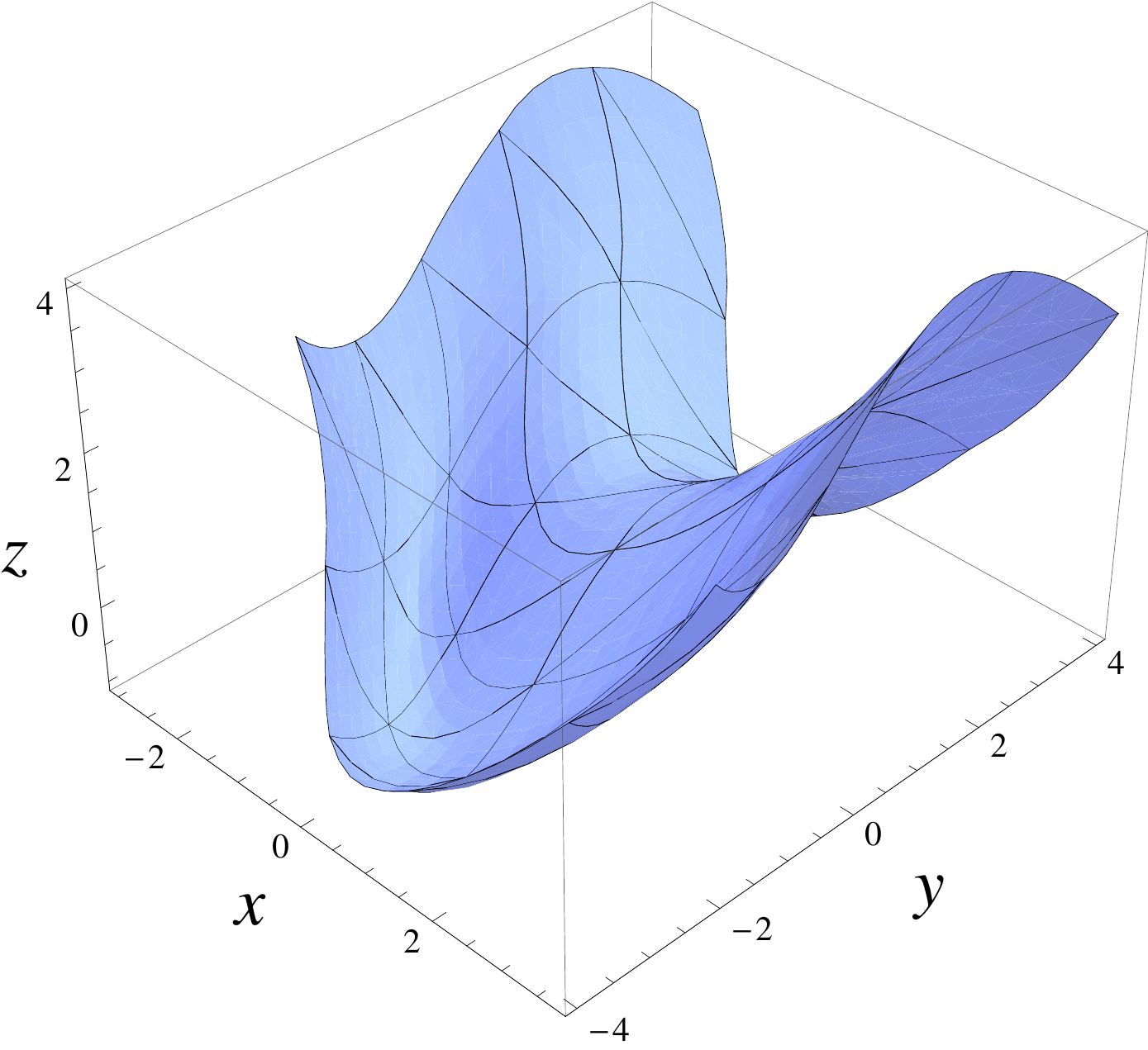}
}
\subfigure[SS scheme, \corrP{$E_H$: 0.14}]{
\includegraphics[height=4.5cm]{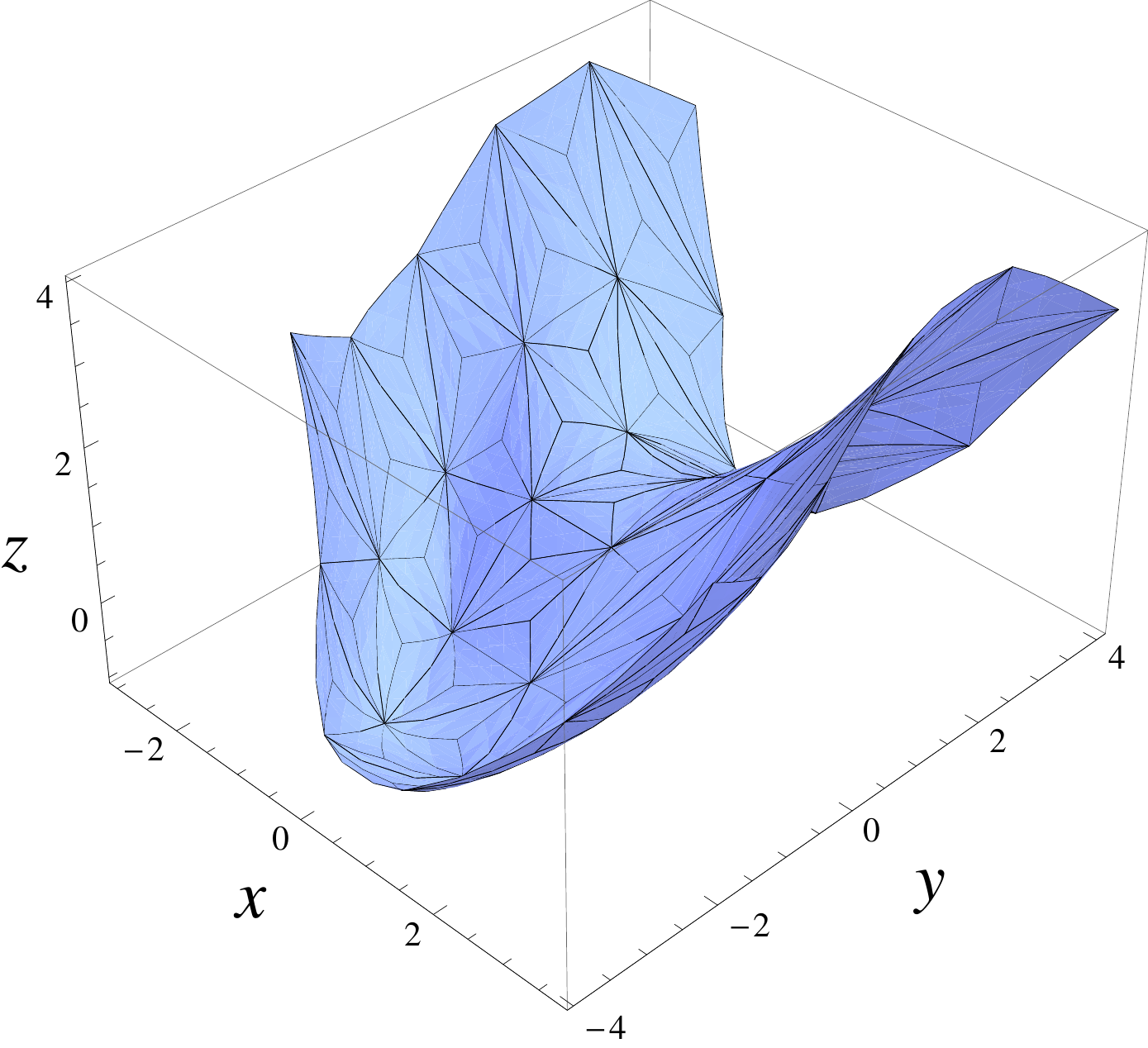}
}
\subfigure[HB scheme, \corrP{$E_H$:  0.36}]{
\includegraphics[height=4.5cm]{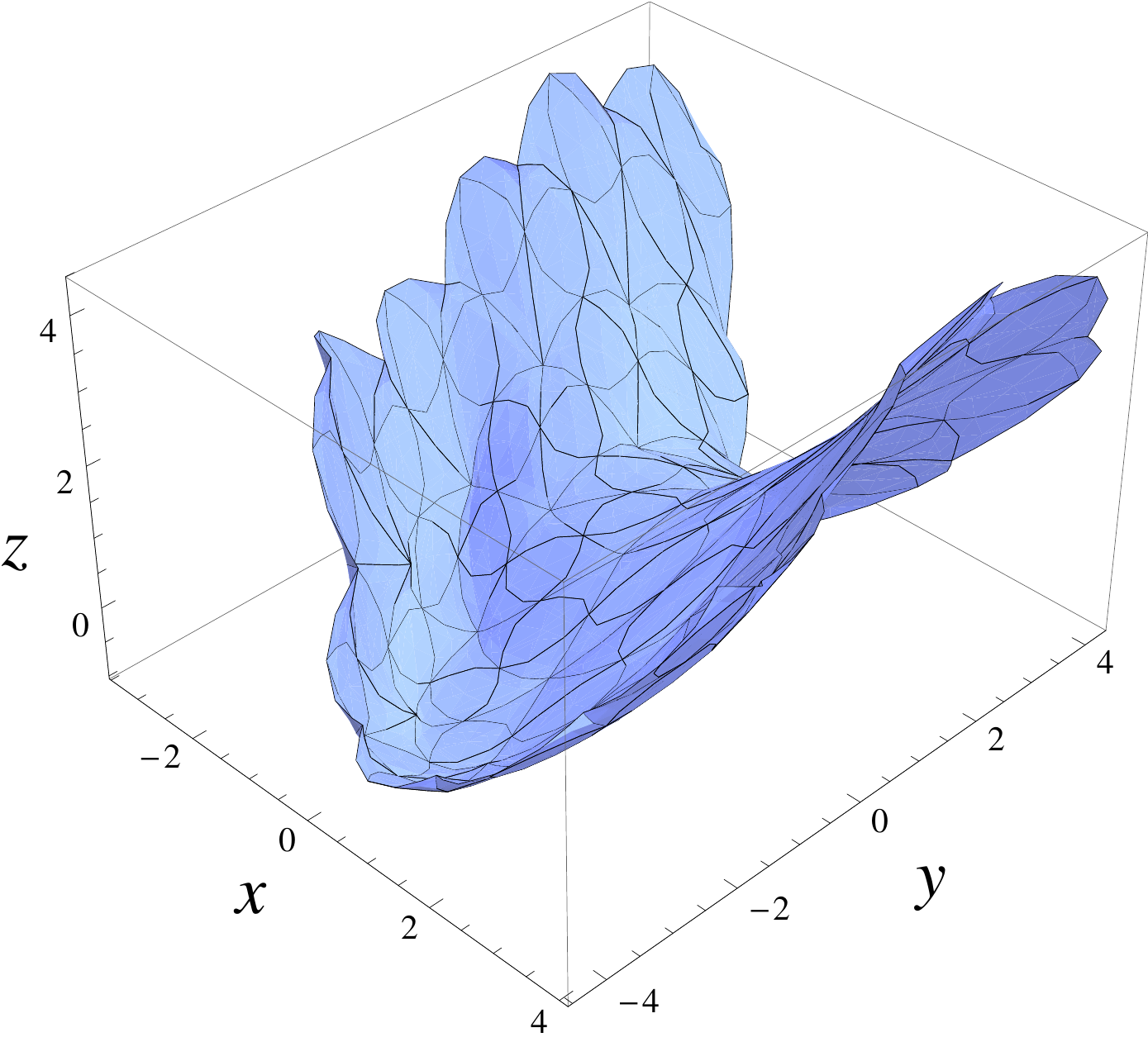}
}
\caption{Approximation of parametric function $\bfm f$ by different interpolation schemes. \corrP{Hausdorff errors ($E_H$) are shown next to the plots.}}
\label{fig:freeFormApprox}
\end{figure}

In the next example, we approximate \corrP{an open} free-form surface defined by a vector function $\bfm f$, 
\begin{align*}
& \bfm f  : [-3,\,3]^2 \to \RR^3,\quad
& \bfm f(u,\, v)  :=   \left(u +  \frac{v^2}{12},\, v - \cos(u),\, \frac 1 3 u^2 + \sin(v) \right).
\end{align*}
Again, we approximate the surface by \corrP{different} interpolation schemes. Plots are shown in Fig.~\ref{fig:freeFormApprox}. As in the first test, our two schemes give smaller Hausdorff distance error than \corrP{SS and HB} methods. 
The referential surface gives a very accurate estimate of the final shape of our two interpolants. 
\corrP{To construct an open surface HB approximant, the spline is constructed on a bigger domain and only relevant patches of the surface are presented.}


\subsection{Approximation of a scalar function}
\begin{figure}[!htb]
\centering
\subfigure[Surface defined by $f$]{
\includegraphics[trim=0cm .0cm 0cm 0cm, clip=true, height=3.75cm]{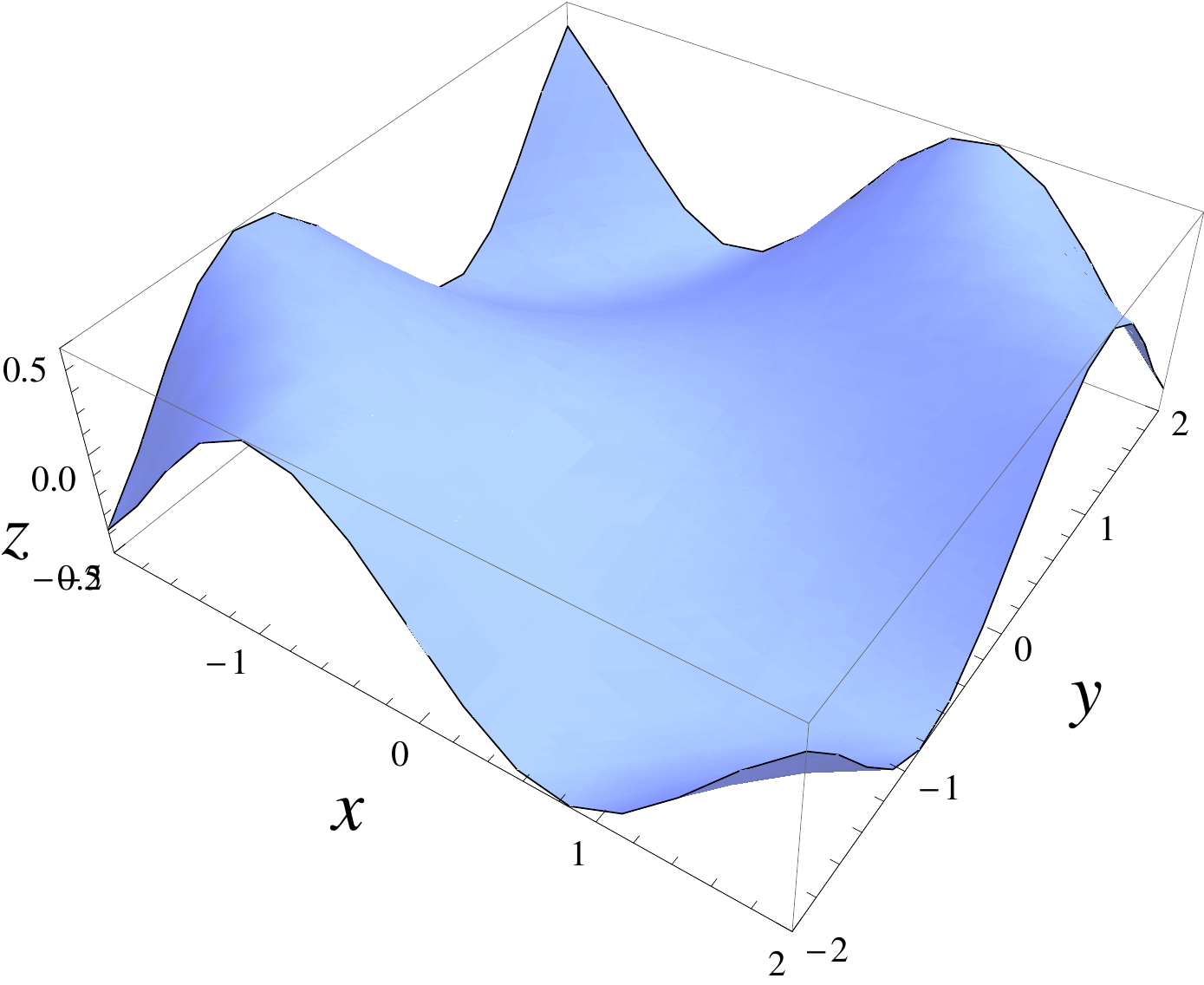}
}
\subfigure[$C^1$ quintic approximant, {$E_{\max}:~0.060$}]{
\includegraphics[trim=-.75cm .0cm -.75cm 0cm, clip=true, height=3.75cm]{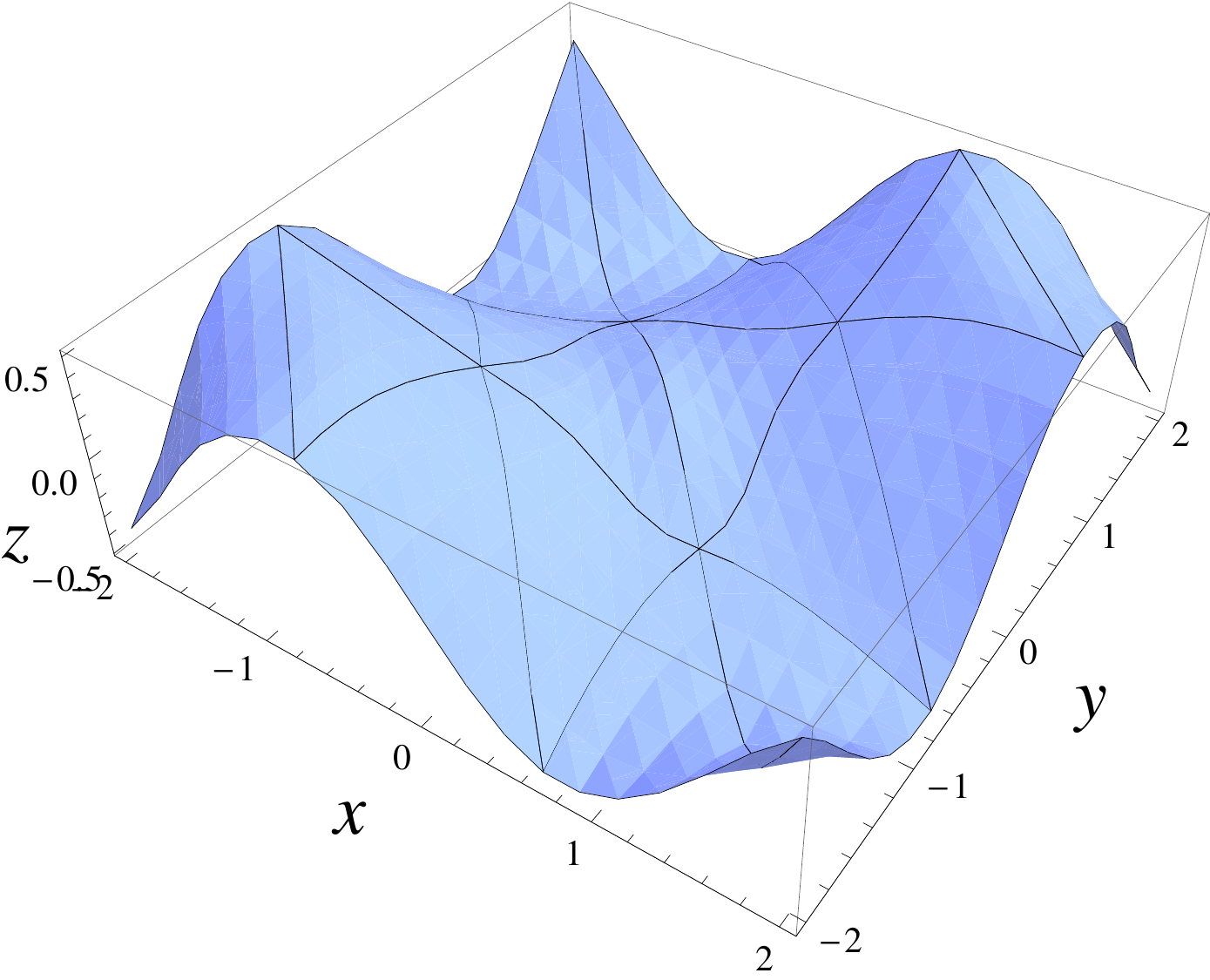}
}
\subfigure[Argyris approximant, {$E_{\max}:~0.015$}]{
\includegraphics[trim=-.75cm .0cm -.75cm 0cm, clip=true, height=3.75cm]{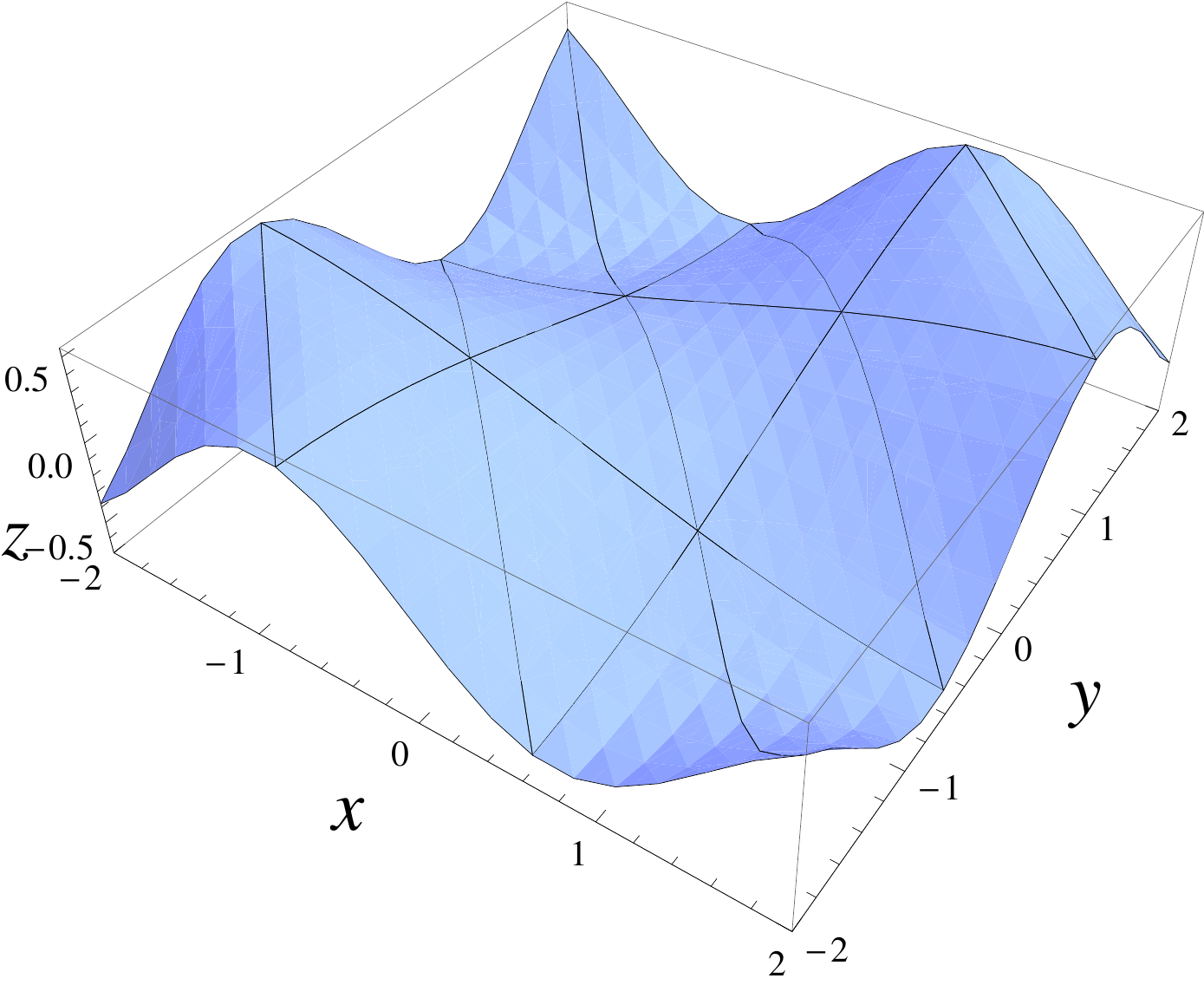}
}
\caption{Approximation of function f by the parametric $C^1$ quintic scheme and the functional Argyris interpolant.
\corrP{Maximal z-errors ($E_{\max}$) are shown next to the plots.}
}
\label{fig:sinApprox}
\end{figure}

In the \corrP{following} example our scheme is compared against functional Argyris interpolant. We approximate a function $f$,
\begin{align*}
&f:[-2,\,2]^2 \to \RR,\quad 
&f(u,\,v) := 1/2\, \sin(u\, v).
\end{align*}
The resulting interpolants are visually almost indistinguishable (Fig.~\ref{fig:sinApprox}). Better accuracy of the Argyris element is expected since it interpolates a much larger set of scalar data that are related to the parameterization, whereas the $C^1$ quintic element interpolates only geometric data.

\subsection{Numerical test of the approximation order}

\begin{figure}[!htb]
\centering
\subfigure[Convergence: sphere approximation]{
\includegraphics[trim=-1cm .0cm -1cm 0cm, height=4.25cm]{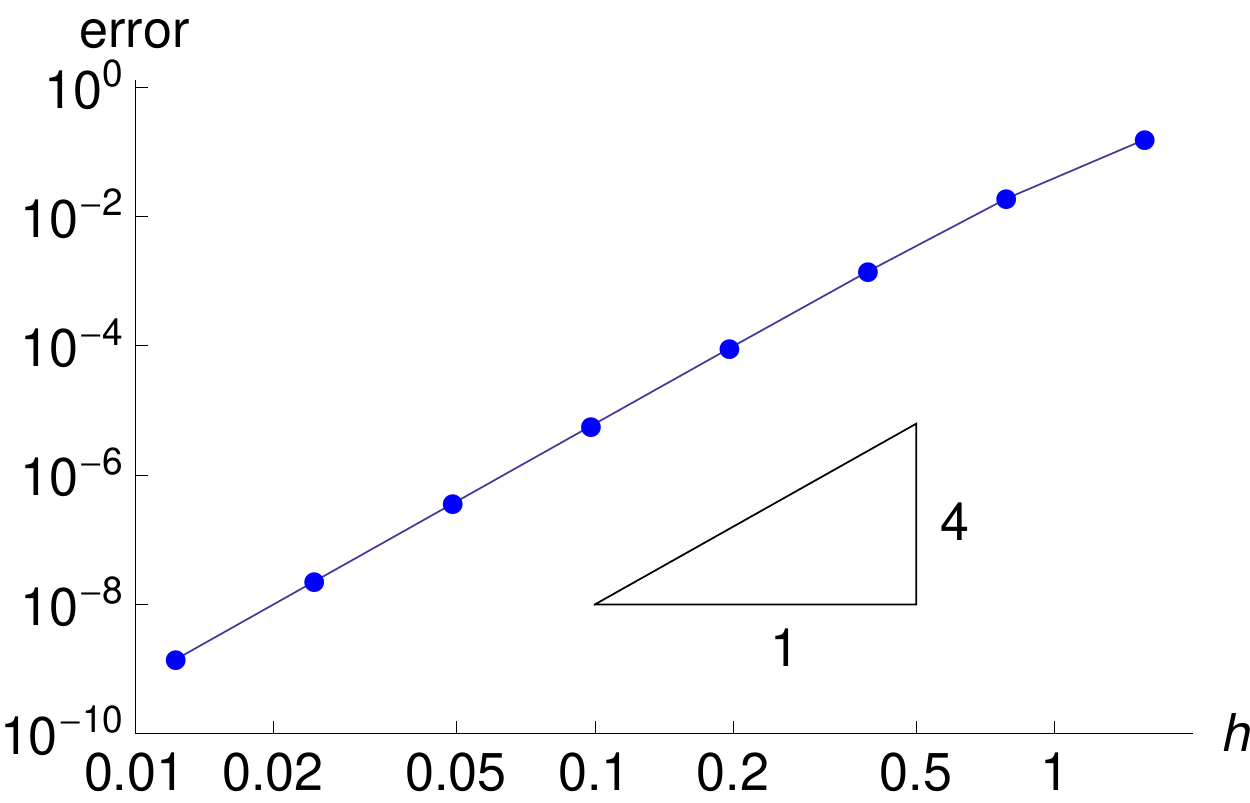}
}
\subfigure[Convergence: approximation of function $f$]{
\includegraphics[trim=-1cm .0cm -1cm 0cm, height=4.25cm]{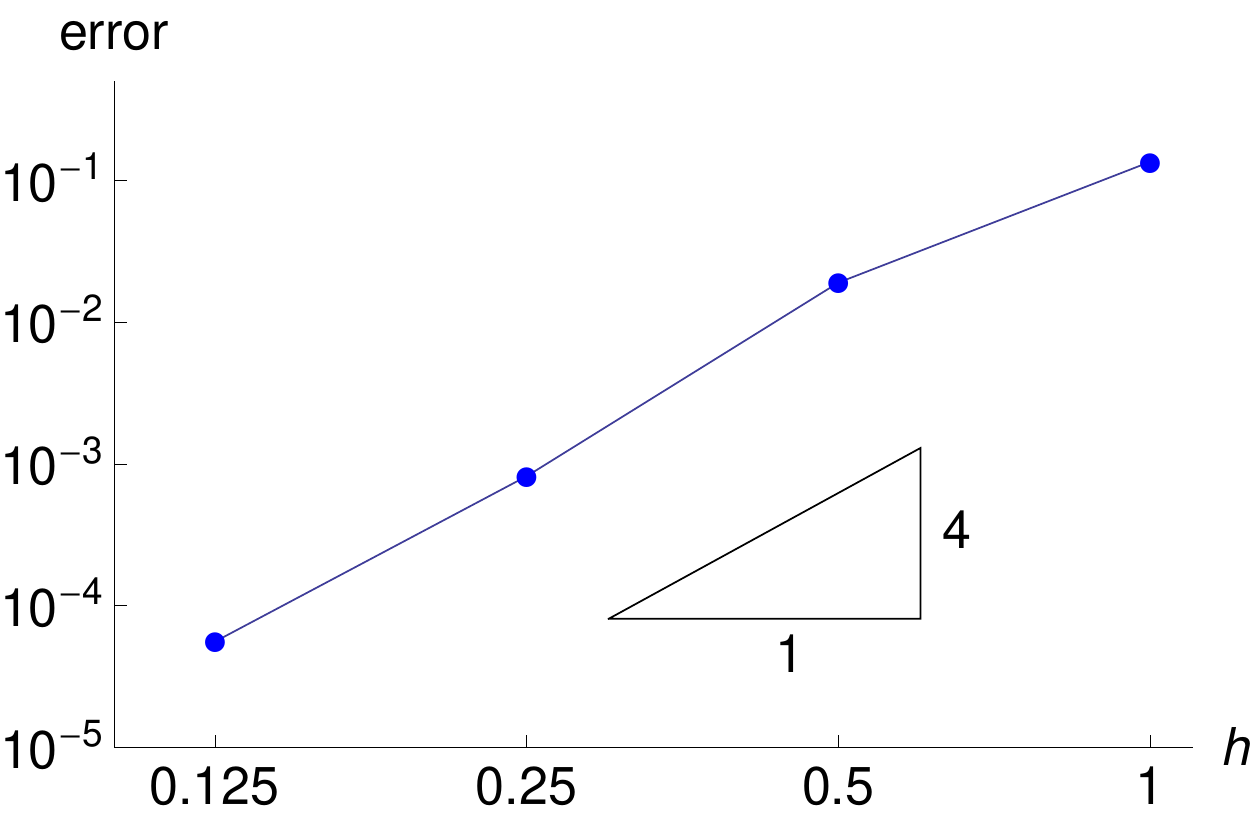}
}
\caption{Convergence plots: error versus the size of triangles ($h$) in triangulation. 
}
\label{fig:convergencePlot}
\end{figure}

\corr{
Since not all of the control points in our scheme are used for approximation, we cannot expect the optimal convergence rate in the general case. The scheme reproduces linear functions, hence the approximation order is at least 2. On each patch 27 scalar data are used for approximation, a number which is close to 30 scalar degrees of freedom of a cubic parametric patch. Therefore, we can speculate that our scheme will approximate well cubic patches and the order of approximation should be 4 in the majority of cases. This observation is confirmed with two numerical tests. In the first one we approximate a unit sphere and in the second function $f=1/2\, \sin(u\, v)$ from the previous example. For the unit sphere case we compare the radial distance error and in the latter case the Hausdorff distance error against the size of triangles in a triangulation. Convergence plots are depicted in Fig.~\ref{fig:convergencePlot}.
}

\section{Conclusions}

In the paper we present a novel Hermite parametric interpolation scheme on triangulations.  \corr{To construct the interpolant, small local systems of equations need to be solved.} 
Two variants of the scheme are derived: $C^1$ quintic and $G^1$ octic. The first has fewer degrees of freedom but needs an explicitly given underlying domain triangulation. The second is of higher polynomial degree but can approximate surfaces of arbitrary topology. The scheme focuses on the geometrical construction of good boundary curves - an important feature to obtain a good approximation surface. $C^2$ smoothness conditions are imposed at the triangle vertices to overcome the twist compatibility problem.

Numerical examples show that our schemes produce approximants with small distance errors and visually satisfying shapes, even when a small number of \corr{patches} is used. On a denser grid of interpolation data, the parametric patches visually resemble the functional patches of Argyris elements.\\

Shape parameters of the interpolant are defined from a referential surface that is based on one step of the Butterfly subdivision scheme. The referential surface gives a basic outline of control points in the space. In the future, the shape parameters could be additionally optimized by applying some energy minimization technique. In that case, control points obtained from the referential surface can be treated as a good starting set of parameters for the optimization procedure. Proper exploitation of the additional six interior control points in $G^1$ octic patches remains an open problem for future work. Combining both scheme variations would results in an adaptive and robust scheme with small number of degrees of freedom. For example, $G^1$ continuity could be applied only around extraordinary vertices. \corr{Another interesting but challenging problem would be to modify parts of the scheme to get the optimal (or near optimal) convergence order while maintaining robustness and desired geometric properties of the scheme.}




\bibliographystyle{elsarticle-num}

\end{document}